# Instability and non-uniqueness in the Cauchy problem for the Euler equations of an ideal incompressible fluid. Part I

Misha Vishik

Department of Mathematics

The University of Texas at Austin

Austin, TX 78712, USA

## § 0. Introduction

In this paper we prove non-uniqueness of the solution to Cauchy problem of the Euler equations of an ideal incompressible fluid in plane with vorticity in some $L^Q(\mathbb{R}^2)$ space. The lack of uniqueness is of the symmetry breaking type, with the radially symmetric external force in $L^1_{loc}([0,\infty); L^Q(\mathbb{R}^2))$.

The weak solution with the initial data in $L^Q(\mathbb{R}^2)$ was defined by V.I. Yudovich in [Y1] and has been constructed by R.DiPerna and A.Majda in [DPM]. In [Y1] the basic uniqueness result for a solution of the initial value problem of the Euler equations of an ideal incompressible fluid in dimension 2 has been proved for the initial vorticity in $L^\infty(\mathbb{R}^2)$. With several later improvements, this result remains the standard reference for the uniqueness problem for the Euler equations in dimension 2.

In his later publications, V.I.Yudovich stressed it would be desirable to extend his uniqueness result to a wider class of the initial data, if possible (see, e.g., [Y3]). In his paper [Y2] he proved uniqueness for the Euler equations in dimension 2 in case the initial vorticity belongs to the intersection of all $L^Q(\mathbb{R}^2), Q \geq Q_0$ $(but\ not\ to\ L^\infty(\mathbb{R}^2))$ with a moderate growth of the $L^Q(\mathbb{R}^2)$ −norm of the initial vorticity. For example, the growth like $O(logQ)$ as $Q \to \infty$ is allowed, but $O(Q)$ in general is not. Therefore, the singularity of the initial vorticity of the type $log|x|$ is not permitted, but $log\,log|x|$ can be handled.

In [V1] the uniqueness result has been proved with vorticity in a borderline Besov space $B_\Gamma$ defined in terms of its Littlewood – Paley decomposition. In terms of admissible singularities this uniqueness theorem allows logarithmic singularities. But the existence theorem has been only proved for the Besov spaces $B_\Gamma$ that do not include this type of singularity. For example, if the vorticity is bounded in time with values in space $BMO$, such a solution is unique.

In V.Scheffer [S] the first remarkable example of non-uniqueness of weak solutions to the Euler equations with velocity with a compact support in $L^2(\mathbb{R}^2 \times \mathbb{R})$ has been presented. A.Shnirelman [Sh1] provided an entirely different and simplified construction of an example of non-uniqueness. In [Sh2] A.Shnirelman constructed an $L^2$ distributional solution with the energy decreasing in time.

In a stunning development C.De Lellis and L.Szekelyhidi [DLS] constructed a solution with a compact support with velocity in $(L^\infty(\mathbb{R}^n \times \mathbb{R}))^n$ and bounded pressure with a compact support. This holds for an arbitrary space dimension n. Their result was based on an entirely different idea of differential inclusions. This direction remains a very active area of research.

In contrast, our solutions are in classes with vorticity in some Lebesgue space $L^Q(\mathbb{R}^2)$, therefore the corresponding velocity is locally Holder in space and time with $Q > 2$. Our approach also describes the resulting non-trivial solution of the Euler equation rather explicitly in terms of self-similar asymptotics up to a small and controlled error. The central idea of the paper, that took years of thinking, is to construct a family of solutions depending on a small parameter $\varepsilon > 0$ to the related system and use the fundamental scaling properties of the Euler equations to pass to the limit as $\varepsilon \to 0$.

In a certain sense, the results of the paper are at variance with the philosophy of V.I.Yudovich explained in his last publications. According to Y.I.Yudovich, the two seemingly unrelated problems, the uniqueness for the Euler equations (i.e., the problem of PDE's) and uniqueness of the Lagrangian flow (the problem of dynamical systems theory), corresponding to a velocity in the certain class of functions,- are deeply related. According to V.I.Yudovich, these two different problems ought to be considered in the same class of functions, and uniqueness holds or does not hold for both problems simultaneously. This philosophy turned out to be true for the classes of uniqueness, discussed in [Y2] and in [V1]. However, the results of this paper demonstrate that this hypothetical correspondence fails for a wider class of the Euler flows with vorticity in Lebesque spaces. Indeed, the well-known remarkable results of P.-L.Lions and R.DiPerna [LDP] (and many subsequent papers on the subject) provide such a uniqueness result for the Lagrangian flow, corresponding to the solenoidal velocity in the Sobolev space $W^{1,Q}(\mathbb{R}^2)$. But the Euler uniqueness fails in this class. Nonetheless, we hope that the philosophy of V.I.Yudovich can still be validated, perhaps under a broader interpretation.

In Part II, the technical result on the linear instability of a certain class of radially symmetric flows is proved (see Theorem 1.1). The problem without an external force is work in progress (Part III of the present paper).

The constant $C$ in the inequalities below denotes a generic constant, the exact value of $C$ may change from line to line.

# §1. Linearized stability problem and the eigenfunctions

Let $G(|x|)$, $x \in \mathbb{R}^2$, $x = (x_1, x_2)$ be radially symmetric $C^\infty(\mathbb{R}^2)$ distribution of vorticity in plane. Define

$$V(x) = \frac{x^\perp}{|x|^2} \int_0^{|x|} s\, G(s)\, ds\,, \qquad\qquad x^\perp = (-x_2, x_1). \tag{1.1}$$

Then $curl\, V(x) = \partial_1 V_2 - \partial_2 V_1 = G(|x|)$.

We choose $\alpha \in (0,2)$ and assume

$$G(s) = s^{-\alpha}, for\ s \geq M > 0\,, \tag{1.2}$$

$M$ being a sufficiently large constant.

Any velocity field (1.1) does (formally) satisfy the Euler equations of an ideal incompressible fluid.

To study the spectrum of small perturbations of this steady-state solution we introduce infinitesimal perturbation

$$w(x) = \partial_x^\perp(e^{im\theta}\,\psi(|x|)), \quad \partial_x^\perp = (-\partial_2, \partial_1)\,; \tag{1.3}$$

$$x_1 = |x|\cos\theta\,, \quad x_2 = |x|\sin\theta\,, \quad \theta \in \mathbb{R}/2\pi\mathbb{Z}\,, \quad x \in \mathbb{R}^2, m \in \mathbb{Z}\backslash\{0\}.$$

Let also

$$e^{im\theta}g(|x|) = curl\, w(x) = \partial_1 w_2 - \partial_2 w_1\,. \tag{1.4}$$

From (1.3), (1.4),

$$e^{im\theta}g(|x|) = \Delta(e^{im\theta}\,\psi(|x|));$$

$$g(s) = \frac{d^2}{ds^2}\psi + \frac{1}{s}\frac{d}{ds}\psi - \frac{m^2}{s^2}\,\psi \tag{1.5}$$

We have,

$$\psi(s) = -\frac{1}{2|m|}\,s^{|m|}\int_s^\infty g(\tau)\tau^{1-|m|}\,d\tau - \frac{1}{2|m|}\,s^{-|m|}\int_0^s g(\tau)\,\tau^{1+|m|}\,d\tau\,. \tag{1.6}$$

Indeed,

$$\frac{d^2}{ds^2}\psi + \frac{1}{s}\frac{d}{ds}\psi - \frac{m^2}{s^2}\psi$$

$$= \frac{1}{s}\frac{d}{ds}s\frac{d}{ds}\psi - \frac{m^2}{s^2}\psi$$

$$= \frac{1}{s}\frac{d}{ds}\left(-\frac{1}{2}s^{|m|}\int_s^\infty g(\tau)\tau^{1-|m|}\,d\tau + \frac{1}{2}s^{-|m|}\int_0^s g(\tau)\,\tau^{1+|m|}\,d\tau\right)$$

$$-\frac{m^2}{s^2}\psi$$

$$= g(s).$$

From (1.4) - (1.6),

$$w(x) = \partial_x^\perp\left(e^{im\theta}\psi(|x|)\right) = \psi'(|x|)e^{im\theta}\frac{x^\perp}{|x|} - im\,\psi(|x|)\,e^{im\theta}\frac{x}{|x|^2} =$$

$$\left(-\frac{|m|}{2|m|}|x|^{|m|-1}\int_{|x|}^\infty g(\tau)\tau^{1-|m|}\,d\tau + \frac{|m|}{2|m|}|x|^{-|m|-1}\int_0^{|x|}g(\tau)\tau^{1+|m|}\,d\tau\right)e^{im\theta}\frac{x^\perp}{|x|}$$

$$-im\,\psi(|x|)\,e^{im\theta}\frac{x}{|x|^2} =$$

$$\left(-\frac{1}{2}|x|^{|m|-2}\int_{|x|}^\infty g(\tau)\tau^{1-|m|}\,d\tau\,e^{im\theta} + \frac{1}{2}|x|^{-|m|-2}\int_0^{|x|}g(\tau)\tau^{1+|m|}\,d\tau\,e^{im\theta}\right)x^\perp$$

$$-im\,\psi(|x|)\,e^{im\theta}\frac{x}{|x|^2}$$

$$= |x|^{|m|-2}\int_{|x|}^\infty g(\tau)\tau^{1-|m|}\,d\tau\,e^{im\theta}\left(-\frac{1}{2}x^\perp + \frac{i}{2}sign\,(m)x\right)$$

$$+ |x|^{-|m|-2}\int_0^{|x|}g(\tau)\tau^{1+|m|}\,d\tau\,e^{im\theta}\left(\frac{1}{2}x^\perp + \frac{i}{2}sign\,(m)x\right)\,. \qquad (1.7)$$

The spectral problem for small oscillations around the steady state solution $V(x)$ looks as follows:

$$\begin{cases}-(V(x),\partial_x)w - (w,\partial_x)V(x) - \partial_x p = \lambda\,w(x)\\ div\,w = 0.\end{cases} \qquad (1.7')$$

We are looking for a nontrivial solution $w \neq 0$ , $w \in L^2(\mathbb{R}^2)$; $Re\,\lambda \neq 0$,

where $w$ is satisfying

$$w(\,R_\theta\,x) = e^{im\theta}\,R_\theta\,w(x), \text{ where } R_\theta = \begin{pmatrix} cos\theta & -sin\theta\\ sin\theta & cos\theta\end{pmatrix},\ \forall\,\theta \in \mathbb{R}/2\pi\mathbb{Z}\,.$$

Then in $S'(\mathbb{R}^2)$ and in $W^{-1,2}(\mathbb{R}^2)$

$$curl\,w\,(R_\theta\,x) = e^{im\theta}curl\,w(x).$$

Therefore,

$$(x^\perp, \partial_x) \, curl \, w = im \, curl \, w$$

in $S'(\mathbb{R}^2)$.

Differentiating the equation for eigenfunctions, we get

$$-(V(x), \partial_x) curl \, w - (w, \partial_x) curl \, V(x) = \lambda \, curl \, w.$$

But $V(x) = \Omega(|x|) x^\perp$, where $\Omega(|x|) \in C^\infty(\mathbb{R}^2) \cap S_{1,0}^{-\alpha}(\mathbb{R}^2)$. Therefore, in $S'(\mathbb{R}^2)$

$$(-im \, \Omega(|x|) - \lambda) \, curl \, w = \, (w, \partial_x) \, G(|x|),$$

where

$$\Omega(s) = = \frac{1}{s^2} \int_0^s \tau \, G(\tau) d\tau \, , \quad s \geq 0. \tag{1.8}$$

But $\{x \to (-im \, \Omega(|x|) - \lambda)\}$ is an invertible multiplier in $S(\mathbb{R}^2)$ and in $S'(\mathbb{R}^2)$. The right side

$$(w, \partial_x) \, G(|x|) \in L^1(\mathbb{R}^2) \cap L^2(\mathbb{R}^2).$$

Therefore,

$$curl \, w(x) = e^{im\theta} g(|x|) \, \in L^1(\mathbb{R}^2) \cap L^2(\mathbb{R}^2).$$

By Sobolev embedding theorem, $w \in L^Q(\mathbb{R}^2), \quad \forall Q \in (2, \infty)$.

Therefore, a unique stream function $e^{im\theta}\psi(|x|)$, such that $w(x) = \partial_x^\perp(e^{im\theta}\psi(|x|)), \psi(0) = 0$,

satisfies $\psi(|x|) \in \, \dot{C}^{1-\frac{2}{Q}}, \, \forall Q \in (2, \infty)$. In particular,

$$|\psi(|x|)| \leq C|x|^{1-\frac{2}{Q}}, \qquad \forall x \in \mathbb{R}^2.$$

In what follows we will have $|m| \geq 2$.

The spectral problem for small oscillations around the steady state solution $V(x)$ in terms of vorticity perturbations looks as follows:

$$-im \, \Omega(|x|) g(|x|) + im \frac{1}{|x|} G'(|x|) \psi(|x|) = \, \lambda \, g(|x|). \tag{1.9}$$

Let

$$E(s) = \frac{1}{s} G'(s). \tag{1.10}$$

From (1.8) - (1.10)

$$\Omega(s) g(s) - E(s)\psi(s) = \mu g(s), \tag{1.11}$$

$$\mu = i\,\frac{\lambda}{m}. \tag{1.12}$$

From (1.2), (1.9), (1.10)

$$E(s) = -\alpha s^{-\alpha-2}, \quad s \geq M; \tag{1.13}$$

$$\Omega(s) = C\alpha\,s^{-2} + \frac{1}{2-\alpha}s^{-\alpha}, \quad s \geq M \tag{1.14}$$

for an appropriate constant $C$.

For the solutions with $Im\,\mu \neq 0$ of the eigenvalue problem for $w(x) \in L^2(\mathbb{R}^2)$,

it follows, that

$$\partial_x\left(e^{im\theta}\psi(|x|)\right) \in L^2(\mathbb{R}^2). \tag{1.15}$$

Therefore,

$$\int_0^\infty s^{-1}\,|\psi(s)|^2 ds\ < \infty, \int_0^\infty s\ |\frac{d\psi}{ds}|^2 ds < \infty,$$

and, consequently, $\psi \in L^\infty(\mathbb{R}^2)$. To see this, we have, for any $s > 0$

$$|\psi(s)|^2 = 2\left|\int_s^\infty s^{-\frac{1}{2}}\,\psi(s)\,s^{\frac{1}{2}}\frac{d\psi}{ds}\,ds\right| \leq C\ < \infty.$$

From (1.15), for any eigenfunction $w(x) \in L^2(\mathbb{R}^2)$ as in (1.7') with $Im\,\mu \neq 0$

$$|g(s)| \leq C(1+s)^{-\alpha-2}, \qquad \forall\,s \geq 0. \tag{1.16}$$

Using (1.6), (1.16) we obtain

$$|g(s)| \leq C(1+s)^{-2\alpha-2}, \qquad \forall\,s \geq 0.$$

Continuing bootstrap, we arrive at the inequality

$$|g(s)| \leq C(1+s)^{-(|m|+2+\alpha)}, \qquad \forall\,s \geq 0. \tag{1.17}$$

Using (1.6), (1.17), (1.11), (1.13), (1.14), we get

$$\mu g(s) = \alpha s^{-\alpha-2}\left(-\frac{1}{2|m|}s^{-|m|}\left[\int_0^\infty g(\tau)\tau^{1+|m|}\,d\tau + O(s^{-\alpha})\right] + O(s^{-|m|-\alpha})\right) + O(s^{-|m|-2-2\alpha})$$

Assuming $\int_0^\infty g(\tau)\tau^{1+|m|}\,d\tau \neq 0$, this implies

$$\mu = -\frac{\alpha}{2|m|}\int_0^\infty g(\tau)\tau^{1+|m|}\,d\tau\ ; \tag{1.18}$$

$$\psi(s) = -\frac{1}{2|m|}s^{-|m|}\int_0^\infty g(\tau)\tau^{1+|m|}\,d\tau + O(s^{-|m|-\alpha})\ as\ s \to \infty, \tag{1.19}$$

provided $g(s)$ is normalized as follows:

$$g(s) \sim s^{-|m|-\alpha-2} \quad as \; s \; \to \infty \; . \tag{1.20}$$

Lemma 1.1. If (1.11) holds with $Im \; \mu \neq 0$, $g(|x|) \in L^2(\mathbb{R}^2)$, and $\int_0^\infty g(\tau)\tau^{1+|m|} \, d\tau = 0$, then $g \equiv 0$.

Proof. Under the assumptions of the Lemma (1.6) implies

$$\psi(s) = -\frac{1}{2|m|} \, s^{|m|} \int_s^\infty g(\tau)\tau^{1-|m|} \; d\tau + \frac{1}{2|m|} \, s^{-|m|} \int_s^\infty g(\tau) \, \tau^{1+|m|} \, d\tau.$$

The same bootstrapping argument as above yields the existence of a constant $C_N$ for any $N \in \mathbb{Z}_+$, such that

$|g(s)| \leq C_N (1+s)^{-N}, \; \forall s \geq 0.$

Since $Im \; \mu \neq 0$, (1.11) implies

$|g(s)| \leq C(1+s)^{-2-\alpha}[s^{|m|} \int_s^\infty |g(\tau)| \; \tau^{|m|+2+\alpha} \; \tau^{-2|m|-1-\alpha} \, d\tau$

$+ s^{-|m|} \int_s^\infty |g(\tau)| \; \tau^{|m|+2+\alpha} \tau^{-1-\alpha} d\tau]. \tag{1.21}$

Define for any $s \geq 0$

$f(s) = \sup_{\tau \geq s} |g(\tau)| \; \tau^{|m|+2+\alpha} \; .$

The inequality (1.21) implies

$$f(s) \leq Cs^{-\alpha}f(s), \qquad \forall s > 0.$$

Thus $g(s) = 0$ for $s \geq C$, where the constant $C$ is sufficiently large. Gronwall's inequality, applied to (1.21) yields the statement of the lemma. QED.

We study the eigenfunction equation (1.11) with $Im \; \mu > 0$. From (1.5), (1.11)

$-\frac{d^2}{ds^2} \psi \; -\frac{1}{s}\frac{d}{ds}\psi + \frac{m^2}{s^2} \psi + \frac{E(s)}{\Omega(s)-\mu}\psi = 0 \; . \tag{1.22}$

For the first observation about (1.22) we introduce a new variable $\; t = \log s$, $t \in \mathbb{R}$.

Then (1.22) becomes

$L\psi \equiv -\frac{d^2}{dt^2} \psi \; + \; m^2\psi + \frac{A(t)}{R(t)-\mu}\psi = 0, \tag{1.23}$

where

$A(t) = \; e^{2t}E(e^t), \tag{1.24}$

$R(t) = \Omega(e^t) \; . \tag{1.25}$

Therefore,

$$A(t) = \frac{d}{dt} G(e^t); \tag{1.26}$$

$$R(t) = \int_{-\infty}^{t} e^{-2(t-\tau)} G(e^\tau) d\tau; \tag{1.27}$$

$$A(t) = R''(t) + 2\,R'(t), \quad (.)' = \frac{d}{dt}\,(.)\;. \tag{1.28}$$

From (1.2), (1.14)

$$A(t) = -\alpha\,e^{-\alpha t}, \qquad t \geq \log M\,; \tag{1.29}$$

$$R(t) = C\alpha\,e^{-2t} + \frac{1}{2-\alpha}e^{-\alpha t}, \quad t \geq \log M\,; \tag{1.30}$$

We will prove the existence of a nontrivial solution to (1.23), (1.28) with

$$Im\,\mu > 0 \tag{1.31}$$

for a certain class of potentials $R(t)$ that we are going to describe now.

Class $\mathcal{C}$.

We denote by $\mathcal{C}$ the following class of functions $R \colon \mathbb{R} \to (0, \infty)$

$(i)\;\; R \in C^\infty(\mathbb{R}),\;\; R'(t) < 0\;\; and$

$$\int_{-\infty}^{t} e^{2\tau} A(\tau) d\tau < 0 \;\; \forall\, t \in \mathbb{R}; \tag{1.32}$$

$(ii)\;\; R(t) = R(-\infty) - c_0 e^{2t}, for\; some\; c_0 > 0, \qquad \forall t \leq \log M_1,\; where\;\; 0 < M_1 < M. \tag{1.33}$

This implies

$$G(s) = 2R(-\infty) - 4c_0\,s^2, \forall\, s\, \in [0, M_1].$$

$(iii)\;\; A(t)\; has\; exactly\; 2\; zeroes, \qquad a < b,\;\; A'(a) > 0, A'(b) < 0.$

Remark 1.  The inequality (1.32) holds if and only if

$$\int_{0}^{e^b} s^2 G'(s) ds\; < 0.$$

Remark 2.  The inequality (1.32) implies

$$R'(t) < 0 \;\; \forall\, t \in \mathbb{R}.$$

Indeed, assuming (1.32), we have from (1.27)

$$R'(t) = G(e^t) - 2\int_{-\infty}^{t} e^{-2(t-\tau)} G(e^\tau) d\tau$$

$$= G(e^t) - e^{-2(t-\tau)} \left. G(e^\tau) \right|_{\tau=-\infty}^{\tau=t} + \int_{-\infty}^{t} e^{-2(t-\tau)} \frac{d}{d\tau} G(e^\tau) d\tau$$

$$= e^{-2t} \int_{-\infty}^{t} e^{3\tau} G'(e^\tau) d\tau$$

$$= e^{-2t} \int_{0}^{e^t} s^2 G'(s) ds \ < 0. \tag{1.34}$$

For $m > 0$ we define the kernel

$$K_m(\xi, \eta) = \frac{1}{2m} e^{-m|\xi - \eta|}, \qquad where \ \xi, \eta \in \mathbb{R}. \tag{1.35}$$

This is the Green's function for the operator

$$-\frac{d^2}{dt^2} + m^2$$

on $\mathbb{R}$.

Indeed, for $f \in L^p(\mathbb{R})$ with some $p \in [1, \infty]$, define

$$\psi(t) = \int K_m(\xi, \eta) f(\eta) \, d\eta$$

$$= e^{mt} \frac{1}{2m} \int_{t}^{\infty} e^{-m\eta} f(\eta) d\eta + e^{-mt} \frac{1}{2m} \int_{-\infty}^{t} e^{m\eta} f(\eta) d\eta.$$

Then,

$$-\psi''(t) + m^2 \psi(t) = -[\frac{1}{2} e^{mt} \int_{t}^{\infty} e^{-m\eta} f(\eta) d\eta - \frac{1}{2m} f(t) - \frac{1}{2} e^{-mt} \int_{-\infty}^{t} e^{m\eta} f(\eta) d\eta$$

$$+ \frac{1}{2m} f(t)]' + m^2 \psi(t)$$

$$= -\frac{m}{2} e^{mt} \int_{t}^{\infty} e^{-m\eta} f(\eta) d\eta + \frac{1}{2} f(t) - \frac{m}{2} e^{-mt} \frac{1}{2} \int_{-\infty}^{t} e^{m\eta} f(\eta) d\eta + \frac{1}{2} f(t) + m^2 \psi(t)$$

$$= f(t),$$

as claimed.

The eigenfunction equation (1.23) with $\psi \in L^p(\mathbb{R}^2)$, $p \in [1, \infty]$ can be written as follows:

$$\psi(t) + \int K_m(t,\xi) \frac{A(\xi)}{\Omega(\xi)-\mu} \psi(\xi) d\xi = 0. \tag{1.36}$$

Notice that $\psi \in W^{1,2}(\mathbb{R}, dt)$ is equivalent to $\partial_x \left( e^{im\theta} \psi(\log|x|) \right) \in L^2(\mathbb{R}^2)$, where $m \neq 0$ is an integer.

Indeed, the radial component of the gradient gives

$$\int\limits_0^\infty s^{-2} |\psi'(\log s)|^2 s \, ds \, < \infty.$$

But this integral is equal to

$$\int\limits_{-\infty}^\infty |\psi'(t)|^2 dt.$$

The tangential component of the gradient gives

$$m^2 \int\limits_0^\infty s^{-1} |\psi(\log s)|^2 \, ds < \infty.$$

After changing variables this yields

$$m^2 \int\limits_{-\infty}^\infty |\psi(t)|^2 \, dt < \infty.$$

We define for every real $m > 0$ the set

$\boldsymbol{\mathcal{U}}_m = \{\mu \in \mathbb{C} | \, Im \,\, \mu > 0 \,\, and \,\, \exists \, \psi = \psi_{m,\mu} \in L^2(\mathbb{R}), \psi \neq 0 \,\, solving \,\, the \,\, equation \,\, (1.36)\}.$

Here is the essential technical result on linear instability of radially symmetric flows.

Theorem 1.1. *For any $\alpha \in (0,2)$ there exists a function $R \in \mathcal{C}$ and an integer $m \geq 2$, so that $\# \,\, \boldsymbol{\mathcal{U}}_m = 1$ and $\boldsymbol{\mathcal{U}}_l = \emptyset$ for any positive integer $l > m$.*

Remark. Any solution $\psi \in L^2(\mathbb{R})$ of the equation (1.36) also belongs to space $W^{1,2}(\mathbb{R})$.

The proof of Theorem 1.1 will be given in part II of the paper.

# § 2. Stability problem for vorticity

Here we study the following two stability problems:

$$\begin{cases} -(V, \partial_x)w - (w, \partial_x)V - \partial_x p = \lambda\, w(x) \\ div\, w = 0; \;\; w = (w_1, w_2), \;\; x \in \mathbb{R}^2, \end{cases} \qquad (2.1)$$

and

$$\begin{cases} -(V, \partial_x)\omega - (w, \partial_x)curl\, V = \lambda\, w(x) \\ div\, w = 0 ; \;\; w = \mathcal{K} * \omega. \end{cases} \qquad (2.2)$$

Here $\mathcal{K}$ denotes the Biot-Savart kernel.

We assume that

$$V(x) = \Omega(|x|)x^{\perp}, \qquad x \in \mathbb{R}^2; \qquad (2.3)$$

$$\Omega(|x|) \in C^{\infty}(\mathbb{R}^2); \qquad (2.4)$$

$$\Omega(s) = Cs^{-2} + \frac{1}{2-\alpha}s^{-\alpha}, \qquad s \geq M > 0, where\; \alpha \in (0,1). \qquad (2.5)$$

For any function space the subscript $m$, $(\, m = 1,2,3, \dots\,)$ refers to a subspace of functions (distributions) invariant under rotation by $\frac{2\pi}{m}$ .

For example, $w \in L_m^Q(\mathbb{R}^2)$ means $w = (w_1, w_2) \in \left(L^Q(\mathbb{R}^2)\right)^2$ and

$$w(x) = R_{\frac{2\pi}{m}} w \left( R_{\frac{2\pi}{m}}^{-1} x \right) \quad \forall x \in \mathbb{R}^2.$$

Likewise, $\omega \in L_m^Q(\mathbb{R}^2)$ means $\omega \in L^Q(\mathbb{R}^2)$ and

$$\omega(x) = \omega \left( R_{\frac{2\pi}{m}}^{-1} x \right) \; \forall x \in \mathbb{R}^2.$$

First, we consider the operator

$$\mathfrak{A}\colon \omega \to (w, \partial_y)\, curl\, V(y); \;\; \omega \in L_m^2(\mathbb{R}^2).$$

This operator is defined as follows. If $\omega \in L_m^{Q_1}(\mathbb{R}^2) \cap L_m^2(\mathbb{R}^2)$ where $Q_1 \in (1,2)$, we define

$$w = \mathcal{K} * \omega$$

and by Hardy-Littlewood-Sobolev inequality

$$\|w\|_{L^{Q_1^*}} \leq C \, \|\omega\|_{L^{Q_1}}, \qquad Q_1^* = \frac{2Q_1}{2 - Q_1} \; . \qquad (2.6)$$

Since $\partial_y \, curl \, V(y) \in L^{Q_1'}$, $\quad Q_1' = \frac{Q_1}{Q_1-1}$ ,

$\|\mathfrak{A}\omega\|_{L^2} \leq C \|w\|_{L^{Q_1^*}} \leq C \|\omega\|_{L^{Q_1}}$ ,

Since $\quad \frac{1}{Q_1^*} + \frac{1}{Q_1'} = 1 - \frac{1}{Q_1} + \frac{1}{Q_1} - \frac{1}{2} = \frac{1}{2}$ .

Therefore, $\mathfrak{A}\omega$ is defined on a dense subspace

$$L^2(\mathbb{R}^2) \cap L^{Q_1}(\mathbb{R}^2) \subset L^2(\mathbb{R}^2).$$

Since $curl \, V$ is radially symmetric,

$$\mathfrak{A}: L^{Q_1}_m(\mathbb{R}^2) \cap L^2_m(\mathbb{R}^2) \to L^2_m(\mathbb{R}^2), \qquad m = 1,2,3,\dots .$$

Proposition 2.1. Let $m \geq 2$. Then $\mathfrak{A}$ extends by continuity to a bounded operator

$$\mathfrak{A}: L^2_m(\mathbb{R}^2) \to L^2_m(\mathbb{R}^2).$$

This operator is compact: $\mathfrak{A} \in \mathfrak{S}_\infty(L^2_m(\mathbb{R}^2))$. In fact, $\mathfrak{A} \in \mathfrak{S}_p(L^2_m(\mathbb{R}^2)) \ \forall p > 2$.

Proof. We have

$$\mathfrak{A} = \bigoplus_{l \in \mathbb{Z}} \mathfrak{A}_l ,$$

where

$$\mathfrak{A}_l: \, g_{ml}(|x|)e^{iml\theta} \mapsto \left(w_{ml}(x)e^{iml\theta}, \partial_x\right)curl \, V(x), \qquad m = 2,3,\dots; l \in \mathbb{Z} , \tag{2.7}$$

where

$w_{ml}(x) = |x|^{|ml|-2} \int_{|x|}^\infty g_{ml}(\tau)\tau^{1-|ml|} \, d\tau \left(-\frac{1}{2}x^\perp + \frac{i}{2}sign \, (ml)x\right)$

$+ |x|^{-|ml|-2} \int_0^{|x|} g_{ml}(\tau)\tau^{1+|ml|} \, d\tau \left(\frac{1}{2}x^\perp + \frac{i}{2}sign \, (ml)x\right)$

for $l \neq 0$;

$$w_0(x) = \frac{x^\perp}{|x|^2} \int_0^{|x|} s \, g_0(s)ds \tag{2.8}$$

for $l = 0$.

Since $curl \, V(x)$ is a radial function,

$$\mathfrak{A}_0 = 0.$$

As for $\mathfrak{A}_l$ with $l \neq 0$, the operator $e^{iml\theta} \circ \mathfrak{A}_l \circ e^{-iml\theta}$ is acting with the kernel

$$Q_l(s,\tau) = \begin{cases} \dfrac{i}{2}\, sign\,(l)\, s^{|ml|-2}\, \tau^{-|ml|}\, sG'(s), & \tau \geq s \\ \dfrac{i}{2}\, sign\,(l)\, s^{-|ml|-2}\, \tau^{|ml|}\, sG'(s), & \tau \leq s \end{cases} \qquad (2.9)$$

Therefore,

$$\mathfrak{A}_l(g_{ml}(|x|)e^{iml\theta})(s,\varphi) = e^{iml\varphi}\int_0^\infty Q_l(s,\tau)\, g_{ml}(\tau)\,\tau\, d\tau.$$

We compute the square of Hilbert-Schmidt norm of $\mathfrak{A}_l$,

$$\iint |Q_l(s,\tau)|^2 \, s ds \;\tau d\tau$$

$$\leq C \iint_{\tau \leq s} s^{-2|ml|-2}\,(1+s)^{-2-2\alpha}\,\tau^{2|ml|} s ds \;\tau d\tau$$

$$+ C \iint_{\tau \geq s} s^{2|ml|-2}\,(1+s)^{-2-2\alpha}\tau^{-2|ml|}\, s ds \;\tau d\tau$$

$$\leq C\,\frac{1}{|ml|}\int_0^\infty s^{-2|ml|-2}\,s\,(1+s)^{-2-2\alpha}s^{2|ml|+2}ds$$

$$+ C\,\frac{1}{|ml|}\int_0^\infty s^{2|ml|-2}s\,(1+s)^{-2-2\alpha}\,s^{2-2|ml|}ds$$

$$\leq C\,\frac{1}{|ml|}\int_0^\infty s\,(1+s)^{-2-2\alpha}ds$$

$$\leq C\,\frac{1}{|ml|}\,.$$

If $\{\sigma_j\}_{j=1}^\infty$ is the sequence of singular numbers of $\mathfrak{A} = \bigoplus_{l\in\mathbb{Z}}\mathfrak{A}_l$, then

$$\left(\sum_{j=1}^\infty {\sigma_j}^p\right)^{\frac{1}{p}} \leq \left(\sum_{l\neq 0}\|\mathfrak{A}_l\|_{HS}^p\right)^{\frac{1}{p}} \leq C\left(\sum_{l\neq 0}|ml|^{-\frac{p}{2}}\right)^{\frac{1}{p}} < \infty\;,$$

QED.

In what follows we assume that an integer $m \geq 2$ is fixed.

Proposition 2.2. Let $\lambda \in \mathbb{C}, Re\,\lambda \neq 0$.

1.   Assume, $w(x) \in L_m^2(\mathbb{R}^2)$, is a solution in the sense of $S'(\mathbb{R}^2)$ of the equation (2.1). Then $\omega = curl\, w(x) \in L_m^2(\mathbb{R}^2)$ and solves (2.2), understood as

$$-(V, \partial_x)\omega - \mathfrak{A}\omega = \lambda\omega. \qquad (2.10)$$

In fact, $\omega \in L_m^1(\mathbb{R}^2) \cap L_m^\infty(\mathbb{R}^2)$ and

$$w = \mathcal{K} * \omega. \qquad (2.11)$$

2. Assume $\omega \in L_m^2(\mathbb{R}^2)$ solves (2.10). Then there exists a unique $w(x) \in L_m^2(\mathbb{R}^2)$.

such that $\omega = curl\ w$ and $w$ solves (2.1).

Proof. According to [V2] the multiplicity of $\lambda$ is finite [since the essential spectral radius in $L_m^2(\mathbb{R}^2)$ of the evolution operator associated with the left side of (2.1) is 1]. Therefore, we may assume that $w$ satisfies

$$w(R_\theta\ y) = e^{iml\theta} R_\theta\ w(y), \quad \forall y \in \mathbb{R}^2.$$

In this functional equation $l \neq 0$, since otherwise $\lambda = 0$ or $w = 0$. Therefore the argument above (following (1.8)) implies statement 1.

Assuming (2.10) with $Re\ \lambda \neq 0$, every nontrivial $l$-component of $\omega$ satisfies the same equation. We fix an integer $l \neq 0$ and $\omega = g(|x|)e^{iml\theta} \neq 0$ solving (2.10). Following the argument above,

$$w \in L_m^2(\mathbb{R}^2), \quad curl\ w = \omega.$$

Such a $w$ is obviously unique. QED.

# §3. Perturbed stability problem

For a small parameter $\kappa > 0$ we consider the following spectral problem

$$L_\kappa\ \omega \equiv -(V, \partial_y)\omega - (w, \partial_x)curl\ V + \kappa((\ y, \partial_y)\omega + \omega)) = \lambda w, \qquad (3.1)$$

where $w = \mathcal{K} * \omega$.

We assume that $V$ satisfies $(2.3) - (2.5)$.

We study the limit $\kappa \to 0 +$.

Proposition 3.1. Assume the spectral problem (2.10) has a total of 2 simple eigenvalues in $L_m^2(\mathbb{R}^2)$

($m \geq 2$ is a fixed integer) with $Re\,\lambda > 0$: $\lambda_0$ and $\overline{\lambda_0}$,

where $\lambda_0 = a + bi$; $\quad a, b > 0$.

Let $d \in (0, a)$.

Then, there is a $\kappa_0 > 0$ such that for every $\kappa \in [0, \kappa_0]$ there are exactly 2 simple eigenvalues of the problem (3.1) in $L_m^2(\mathbb{R}^2)$ in half-plane $\{\lambda \in \mathbb{C}|\ Re\,\lambda > d\}$: $\lambda_\kappa$ and $\overline{\lambda_0}$; $Im\,\lambda_\kappa > 0$.

In addition to this $\lim\limits_{\kappa \to 0+} \lambda_\kappa = \lambda_0$.

Proof. Let

$M_\kappa \omega \equiv -(V, \partial_y)\omega + \kappa\big((y, \partial_y)\omega + \omega)\big).$

Given the assumptions on $V$ the operator $M_\kappa$ is skew -Hermitian in $L_m^2(\mathbb{R}^2)$ . Therefore,

$$\|(\lambda - M_\kappa)^{-1}\|_{\mathcal{L}\left(L_m^2(\mathbb{R}^2)\right)} \leq |Re\,\lambda|^{-1},\ Re\,\lambda \neq 0.$$

Assume, $(\omega, \lambda), \omega \neq 0$ is a solution to the eigenvalue problem (3.1), $\omega \in L_m^2(\mathbb{R}^2)$ . Then,

$\omega + R_\lambda(M_\kappa)\,\mathfrak{A}\,\omega = 0,$ \hfill (3.2)

Where $R_\lambda(M_\kappa) = (\lambda - M_\kappa)^{-1} \in \mathcal{L}\big(L_m^2(\mathbb{R}^2)\big),\ \forall \lambda,\ Re\,\lambda \neq 0.$

Therefore, since $\mathfrak{A}$ is bounded in $L_m^2(\mathbb{R}^2)$, there exists a constant $D \in (d, \infty)$ such that for a $\lambda \in \mathbb{C}$, an eigenvalue in (3.1) necessarily $Re\lambda < D$.

Let $\lambda = s + it, \kappa_0 > 0$ be fixed. First, we show that

$$\lim\limits_{|t| \to \infty} \sup\limits_{s \in [d,D], \kappa \in [0,\kappa_0]} \|R_\lambda(M_\kappa)\,\mathfrak{A}\|_{\mathcal{L}\left(L_m^2(\mathbb{R}^2)\right)} = 0.$$ \hfill (3.3)

Let $\varepsilon \in (0,1)$ be given. We have the following inequalities from the Laplace transform formula:

$$R_\lambda(M_\kappa) = \int_0^\infty \exp(-\lambda\tau + M_\kappa\tau)\,d\tau$$

$$= \int_0^T \exp(-\lambda\tau + M_\kappa\tau)\,d\tau + \int_T^\infty \exp(-\lambda\tau + M_\kappa\tau)\,d\tau\,,$$ \hfill (3.4)

$$\|\int_T^\infty \exp(-\lambda\tau + M_\kappa\tau)d\tau\|_{\mathcal{L}\left(L_m^2(\mathbb{R}^2)\right)} \leq \frac{1}{d}e^{-dT}\,,$$ \hfill (3.5)

if $s = Re\lambda \geq d$.

We choose $T$ so that $\exp(-dT) < \varepsilon$.

Since $\mathfrak{A} \in \mathfrak{S}_\infty\left(L^2_m(\mathbb{R}^2)\right)$, there exists an approximation $\mathfrak{A}^{(\varepsilon)} \in \mathcal{L}\left(L^2_m(\mathbb{R}^2)\right)$, $\;rk\,\mathfrak{A}^{(\varepsilon)} < \infty$, such that

$$\|\mathfrak{A} - \mathfrak{A}^{(\varepsilon)}\|_{\mathcal{L}\left(L^2_m(\mathbb{R}^2)\right)} < \varepsilon. \tag{3.6}$$

We may represent $\mathfrak{A}^{(\varepsilon)}$ as follows:

$\mathfrak{A}^{(\varepsilon)} = \sum_{j=1}^{N}(.\,,g_j)f_j \quad$ ,

where $f_j, g_j \in L^2_m(\mathbb{R}^2)$ for $j = 1,2,3,\ldots,N$.

By approximating $g_j, f_j$ in $L^2_m(\mathbb{R}^2)$, we may assume that $g_j, f_j \in C_0^\infty(\mathbb{R}^2), j = 1,2,3,\ldots,N$.

We have:

$$\|R_\lambda(M_\kappa)\mathfrak{A} - R_\lambda(M_\kappa)\mathfrak{A}^{(\varepsilon)}\|_{\mathcal{L}\left(L^2_m(\mathbb{R}^2)\right)} < \varepsilon\, d^{-1}, \qquad Re\,\lambda \geq d. \tag{3.7}$$

Using (3.4) --(3.7) for any $f \in L^2_m(\mathbb{R}^2)$,

$$\|R_\lambda(M_\kappa)\mathfrak{A}f - \sum_{j=1}^{N}\int_0^T \exp(-\lambda\tau + M_\kappa\tau)(f,g_j)f_j\,dt\|_{\mathcal{L}\left(L^2_m(\mathbb{R}^2)\right)}$$

$$\leq \frac{\varepsilon}{d}(\|\mathfrak{A}\|_{\mathcal{L}\left(L^2_m(\mathbb{R}^2)\right)} + 1)\|f\|_{L^2_m(\mathbb{R}^2)} \tag{3.8}$$

Integrating by parts in the left side of (3.8), we get

$$\int_0^T \exp(-\lambda t + M_\kappa t)f_j\,dt = -\frac{1}{\lambda}\exp(-\lambda T)\exp(M_\kappa T)f_j + \frac{1}{\lambda}f_j$$

$+\frac{1}{\lambda}\int_0^T \exp(-\lambda t + M_\kappa t)M_\kappa f_j\,dt.$

Since $M_\kappa$ is skew -Hermitian in $L^2_m(\mathbb{R}^2)$, $Re\,\lambda \geq d$, we get

$$\|\int_0^T \exp(-\lambda t + M_\kappa t)f_j\,dt\|_{L^2_m(\mathbb{R}^2)} \leq 2\,\frac{1}{|Im\,\lambda|}\|f_j\|_{L^2_m(\mathbb{R}^2)} + \frac{1}{d\,|Im\,\lambda|}\|M_\kappa f_j\|_{L^2_m(\mathbb{R}^2)}\,. \tag{3.9}$$

Therefore, from (3.8), (3.9)

$$\|R_\lambda(M_\kappa)\mathfrak{A}\|_{\mathcal{L}\left(L^2_m(\mathbb{R}^2)\right)} \leq \frac{\varepsilon}{d}(\|\mathfrak{A}\|_{\mathcal{L}\left(L^2_m(\mathbb{R}^2)\right)} + 1)$$

$+\frac{1}{|Im\,\lambda|}\sum_{j=1}^{N}\|g_j\|_{L^2_m(\mathbb{R}^2)}(2\|f_j\|_{L^2_m(\mathbb{R}^2)} + \frac{1}{d}\|M_\kappa f_j\|_{L^2_m(\mathbb{R}^2)}), \quad Re\,\lambda \geq d.$

Since $\kappa \in [0, \kappa_0]$, $(V, \partial_y) f_j$, $(y, \partial_y) f_j \in L^2_m(\mathbb{R}^2)$, we may choose $|Im \lambda|$ large enough, so that this expression is less than

$$2 \frac{\varepsilon}{d} \left( \|\mathfrak{A}\|_{\mathcal{L}\left(L^2_m(\mathbb{R}^2)\right)} + 1 \right).$$

Let $\Gamma$ be any closed piecewise $C^1$ contour, $\Gamma \subset \{z \in \mathbb{C} | Re z \geq d\}$, avoiding the spectrum of $L_0$ in $L^2_m(\mathbb{R}^2)$, i.e. $\lambda_0 \ \overline{\lambda_0} \notin \Gamma$.

The operator

$$1 - R_\lambda(M_0) \, \mathfrak{A} \in \mathcal{L}\left(L^2_m(\mathbb{R}^2)\right)$$

is invertible if $Re \, \lambda > 0$. Indeed, using Proposition 2.1 and assume the contrary. From the Fredholm alternative, there is $\omega \in L^2_m(\mathbb{R}^2)$, so that $\omega \neq 0$ and

$$\omega - R_\lambda(M_0)\mathfrak{A}\omega = 0$$

Therefore, $\omega \in D(M_0) \subset L^2_m(\mathbb{R}^2)$ and

$$\lambda\omega - M_0\omega - \mathfrak{A}\omega = 0.$$

But $M_0 + \mathfrak{A} = L_0$, so $\lambda = \lambda_0$ or $\lambda = \overline{\lambda_0}$, contradicting the assumptions.

We claim that

$$\lim_{\kappa \to 0} \sup_{\lambda \in \Gamma} \| \left( R_\lambda(M_\kappa) - R_\lambda(M_0)\right)\mathfrak{A}\|_{\mathcal{L}\left(L^2_m(\mathbb{R}^2)\right)} = 0. \tag{3.10}$$

The proof of (3.10) is similar to the one above. Given $\varepsilon \in (0,1)$ let $\mathfrak{A}^{(\varepsilon)} \in \mathcal{L}\left(L^2_m(\mathbb{R}^2)\right)$, $rk \, \mathfrak{A}^{(\varepsilon)} < \infty$ be chosen so that (3.6) holds. We may represent $\mathfrak{A}^{(\varepsilon)}$ as follows:

$$\mathfrak{A}^{(\varepsilon)} = \sum_{j=1}^N (., g_j)f_j \quad ,$$

where $f_j, g_j \in L^2_m(\mathbb{R}^2)$ for $j = 1,2,3,\dots,N$. It follows that

$$\|\left(R_\lambda(M_\kappa) - R_\lambda(M_0)\right)\mathfrak{A} - \left(R_\lambda(M_\kappa) - R_\lambda(M_0)\right)\mathfrak{A}^{(\varepsilon)}\|_{\mathcal{L}\left(L^2_m(\mathbb{R}^2)\right)} < 2\varepsilon d^{-1}, Re \, \lambda \geq d, \kappa \in [0, \kappa_0] \tag{3.11}$$

Also,

$$\|R_\lambda(M_\kappa) - \int_0^T \exp(-\lambda t + M_\kappa t)dt\|_{\mathcal{L}\left(L^2_m(\mathbb{R}^2)\right)} \leq d^{-1}e^{-dT} < d^{-1}\varepsilon, \tag{3.12}$$

if T is sufficiently large, so that $e^{-dT} < \varepsilon$.

Therefore,

$$\|\left(R_\lambda(M_\kappa) - R_\lambda(M_0)\right)\mathfrak{A} - \int_0^T e^{-\lambda t}(\exp\{M_\kappa t\} - \exp\{M_0 t\})dt\ \mathfrak{A}^{(\varepsilon)}\|_{\mathcal{L}\left(L_m^2(\mathbb{R}^2)\right)} \leq$$

$$\leq 2d^{-1}\varepsilon + 2d^{-1}\varepsilon\left(\|\mathfrak{A}\|_{\mathcal{L}\left(L_m^2(\mathbb{R}^2)\right)} + 1\right). \tag{3.13}$$

For $f \in L_m^2(\mathbb{R}^2), \exp\{M_\kappa t\}f$ is the unique solution of the initial value problem

$$\partial_t \omega(y,t) = M_\kappa \omega, \ \ \omega(y,0) = f(y). \tag{3.14}$$

Let $g_\kappa^t y$ be the diffeomorphism corresponding to the dynamical system

$$\dot{x} = V(x) - \kappa x, \ \ x(0) = y\,. \tag{3.15}$$

Then

$g_\kappa^t y = x(t), \ \forall y \in \mathbb{R}^2, t \in \mathbb{R}.$

From (3.14), (3.15)

$$(\exp\{M_\kappa t\}f)(x) = e^{\kappa t}f(g_\kappa^{-t}x). \tag{3.16}$$

Therefore,

$$\int_0^T e^{-\lambda t}(\exp\{M_\kappa t\} - \exp\{M_0 t\})\,\mathfrak{A}^{(\varepsilon)}f\,dt$$

$$= \sum_{j=1}^N (f, g_j) \int_0^T \left(f_j(g_\kappa^{-t}x) - f_j(g_0^{-t}x)\right) e^{-\lambda t + \kappa t}dt.$$

We select a $\kappa_0 > 0$ so that for any $\kappa \in [0, \kappa_0]$

$$\|e^{\kappa t}(f_j(g_\kappa^{-t}x) - f_j(g_0^{-t}x))\|_{L_m^2(\mathbb{R}^2)}\|g_j\|_{L_m^2(\mathbb{R}^2)} < N^{-1}\varepsilon, \ \ j = 1,2,3, \dots N.$$

Then,

$$\|\int_0^T e^{-\lambda t}(\exp\{M_\kappa t\} - \exp\{M_0 t\})\,\mathfrak{A}^{(\varepsilon)}\,dt\|_{\mathcal{L}\left(L_m^2(\mathbb{R}^2)\right)} < d^{-1}\varepsilon.$$

From (3.13)

$$\|\left(R_\lambda(M_\kappa) - R_\lambda(M_0)\right)\mathfrak{A}\|_{\mathcal{L}\left(L_m^2(\mathbb{R}^2)\right)} \leq$$

$$\leq 3d^{-1}\varepsilon + 2d^{-1}\varepsilon\left(\|\mathfrak{A}\|_{\mathcal{L}\left(L_m^2(\mathbb{R}^2)\right)} + 1\right),$$

for $Re\ \lambda \geq d$ and $\kappa \in [0, \kappa_0].$

This proves (3.10).

A similar proof yields

$$\lim_{\kappa \to 0} \sup_{\lambda \in \Gamma} \| \mathfrak{A}\left(R_\lambda(M_\kappa) - R_\lambda(M_0)\right)\|_{\mathcal{L}\left(L_m^2(\mathbb{R}^2)\right)} = 0. \tag{3.17}$$

Using (3.3) we find a $T > 0$ and $\kappa_1 > 0$ so that

$$\sup_{s \in [d,D], \kappa \in [0,\kappa_1]} \| R_\lambda(M_\kappa)\mathfrak{A}\|_{\mathcal{L}\left(L_m^2(\mathbb{R}^2)\right)} < \frac{1}{2} \ , \qquad |t| \geq T.$$

Therefore, for $\kappa \in [0,\kappa_1]$

$$\sigma_{L_m^2(\mathbb{R}^2)}(L_\kappa) \cap \{ \lambda \in \mathbb{C} | \ Re \ \lambda \in [d,D], |Im \ \lambda| \geq T\} = \emptyset.$$

Also, trivially,

$$\sigma_{L_m^2(\mathbb{R}^2)}(L_\kappa) \cap \{ \lambda \in \mathbb{C} | \ Re \ \lambda \geq D\} = \emptyset.$$

We construct a closed contour $\Gamma$ with positive orientation as follows:

$$\Gamma = \{s - iT | s \in [d,D]\} \cup \{D + it | t \in [-T,T]\} \cup \{s + iT | s \in [d,D]\} \cup \{d + it | t \in [-T,T]\}.$$

The contour $\Gamma$ contains $\lambda_0$, $\overline{\lambda_0}$ as interior points.

The operator $1 - R_\lambda(M_0)\mathfrak{A} \in \mathcal{L}\left(L_m^2(\mathbb{R}^2)\right)$ is invertible and analytic on $\Gamma$. From (3.10) for any $\delta > 0$ there exists a $\kappa_0 > 0$, such that for $\kappa \in [0,\kappa_0]$ the operator

$$1 - R_\lambda(M_\kappa)\mathfrak{A} \in \mathcal{L}\left(L_m^2(\mathbb{R}^2)\right)$$

is also invertible and

$$\sup_{\lambda \in \Gamma} \| (1 - R_\lambda(M_\kappa)\mathfrak{A})^{-1} - (1 - R_\lambda(M_0)\mathfrak{A})^{-1}\|_{\mathcal{L}\left(L_m^2(\mathbb{R}^2)\right)} < \delta, \quad \forall \, \kappa \in [0,\kappa_0], \ \forall \lambda \in \Gamma. \tag{3.18}$$

Hence, for $\kappa \in [0,\kappa_0], \ \lambda \in \Gamma$

$$R_\lambda(L_\kappa) = (1 - R_\lambda(M_\kappa)\mathfrak{A})^{-1}R_\lambda(M_\kappa) \ \in \ \mathcal{L}\left(L_m^2(\mathbb{R}^2)\right)$$

exists and is analytic on $\Gamma$.

Let for $f \in L_m^2(\mathbb{R}^2), \ \lambda \in \Gamma$,

$$\omega = R_\lambda(L_\kappa)f.$$

Then

$$\omega - R_\lambda(M_\kappa)\mathfrak{A}\omega = R_\lambda(M_\kappa)f.$$

Let $\omega = R_\lambda(M_\kappa)f + \omega_1$.

Then

$$\omega_1 - R_\lambda(M_\kappa)\mathfrak{A}\omega_1 = R_\lambda(M_\kappa)\mathfrak{A}R_\lambda(M_\kappa)f.$$

Since

$$1 - R_\lambda(M_\kappa)\mathfrak{A} \in \mathcal{L}\big(L_m^2(\mathbb{R}^2)\big)$$

is invertible, we get

$$R_\lambda(L_\kappa) = R_\lambda(M_\kappa) + (1 - R_\lambda(M_\kappa)\mathfrak{A})^{-1}R_\lambda(M_\kappa)\mathfrak{A}R_\lambda(M_\kappa), \ \forall \kappa \in [0, \kappa_0], \ \forall \lambda \in \Gamma.$$

We use (3.10), (3.17) -(3.19) to get

$$\big\|\big(R_\lambda(L_\kappa) - R_\lambda(M_\kappa)\big) - \big(R_\lambda(L_0) - R_\lambda(M_0)\big)\big\|_{\mathcal{L}\big(L_m^2(\mathbb{R}^2)\big)} < C\delta, \forall \kappa \in [0, \kappa_0], \ \forall \lambda \in \Gamma.$$

After integrating the resolvent along the contour $\Gamma$

$$\|\tfrac{1}{2\pi i}\oint_\Gamma \ \big(R_\lambda(L_\kappa) - R_\lambda(L_0)\big) \ d\lambda\|_{\mathcal{L}\big(L_m^2(\mathbb{R}^2)\big)} < C\delta|\Gamma|.$$

Indeed,

$$\tfrac{1}{2\pi i}\oint_\Gamma \ R_\lambda(M_\kappa) \ d\lambda = 0 \ , \ \forall \kappa \in \mathbb{R}.$$

This implies there are exactly two eigenvalues of $L_\kappa$ (counted with multiplicities) inside the contour $\Gamma$.

Now let $\Gamma$ be a circle of an arbitrary small radius $r > 0$ centered at $\lambda_0$ or at $\overline{\lambda_0}$ .The same argument proves there is exactly one eigenvalue of the operator $L_\kappa$ inside $\Gamma$ for sufficiently small $\kappa > 0$. This completes the proof of Proposition 3.1. QED.

# §4. Weak solutions of the Euler Equations in dimension 2.

## Vorticity formulation.

We will follow the definition of a weak solution that is restricted. This definition will imply most of the other definitions.

We study the following system of equations:

$$\begin{cases} \partial_t \omega(x,t) = -(v(x,t), \partial_x)\omega(x,t) + Z(x,t); \ x \in \mathbb{R}^2, \ t \in [0,\infty), \\ \quad \omega(x,t) = curl_x \ v(x,t); \ div_x v(x,t) = 0, \\ \quad \quad \omega(x,0) = \omega_0(x) \in L^{Q_0}(\mathbb{R}^2) \cap L^{Q_1}(\mathbb{R}^2). \end{cases} \quad (4.1)$$

We assume

$$\omega(x,t) \in L^\infty_{loc}\big([0,\infty)\big); \ L^{Q_0}(\mathbb{R}^2) \cap L^{Q_1}(\mathbb{R}^2); \quad (4.2)$$

$$Z(x,t) \in L^1_{loc}\big([0,\infty)\big); \ L^{Q_0}(\mathbb{R}^2) \cap L^{Q_1}(\mathbb{R}^2)), \quad (4.3)$$

where

$$Q_0 \in (1,2) \, , Q_1 \in (2,\infty). \quad (4.4)$$

We define

$$v(x,t) = \big(\mathcal{K} *_x \omega(t)\big)(x). \quad (4.5)$$

Lemma 4.1. Convolution with $\mathcal{K}$ is a bounded operator from $L^{Q_0}(\mathbb{R}^2) \cap L^{Q_1}(\mathbb{R}^2)$ into $L^\infty(\mathbb{R}^2)$:

$$\mathcal{K} *_x \ . \ : L^{Q_0}(\mathbb{R}^2) \cap L^{Q_1}(\mathbb{R}^2) \to L^\infty(\mathbb{R}^2).$$

Proof. This is well known. If $\chi \in C_0^\infty(\mathbb{R}^2)$ is such that $\chi|_{B_0(1)} \equiv 1$, we split the integral as follows:

$$\mathcal{K} * \omega = (\chi\mathcal{K}) * \omega + \big((1-\chi)\mathcal{K}\big) * \omega.$$

This representation yields

$$\|\mathcal{K} * \omega\|_{L^\infty(\mathbb{R}^2)} \le \|\chi\mathcal{K}\|_{L^{Q_1'}(\mathbb{R}^2)} \|\omega\|_{L^{Q_1}(\mathbb{R}^2)} + \|(1-\chi)\mathcal{K}\|_{L^{Q_0'}(\mathbb{R}^2)} \|\omega\|_{L^{Q_0}(\mathbb{R}^2)}$$

$$\le C \ \|\omega\|_{L^{Q_1}(\mathbb{R}^2) \cap L^{Q_0}(\mathbb{R}^2)}.$$

Here $Q_1' = \frac{Q_1}{Q_1 - 1}$ , $Q_0' = \frac{Q_0}{Q_0 - 1}$ .

QED.

Remark. (4.5) implies that $div_x v(x,t) = 0$ in $S'(\mathbb{R}^2)$ and $\omega(x,t) = curl_x v(x,t)$ in $S'(\mathbb{R}^2)$.

Let $\psi(x,t)$ be a test function

$$\psi(x,t) \in C_0^\infty\big([0,\infty); S(\mathbb{R}^2)\big). \qquad (4.6)$$

Let

$$\omega_0(x) \in L^{Q_0}(\mathbb{R}^2) \cap L^{Q_1}(\mathbb{R}^2).$$

Definition of a weak solution. A weak solution to the problem (4.1) is a function $\omega\,(x,t)$,

$$\omega(x,t) \in L_{loc}^\infty\big([0,\infty)\big);\ L^{Q_0}(\mathbb{R}^2) \cap L^{Q_1}(\mathbb{R}^2),$$

such that for every test function $\psi(x,t)$ as in (4.6) the following identity holds true:

$$\int \omega_0(x)\,\psi(x,0)\,dx + \int\limits_0^\infty \int \omega(x,t)\,\partial_t\,\psi(x,t)\,dx\,dt$$

$$\int_0^\infty \int \omega(x,t)\,(v(x,t),\partial_x)\psi(x,t)\,dx\,dt + \int_0^\infty \int Z(x,t)\,\psi(x,t)\,dx\,dt = 0. \qquad (4.7)$$

The integrals in the left side of (4.7) are meaningful because of (4.2)-(4.6) and Lemma 4.1.

We will consider a *variant* of this definition for each $m = 1,2,3,\dots$ .

In this context

$$\omega_0(x) \in L_m^{Q_0}(\mathbb{R}^2) \cap L_m^{Q_1}(\mathbb{R}^2),$$

$$\omega(x,t) \in L_{loc}^\infty\big([0,\infty)\big); L_m^{Q_0}(\mathbb{R}^2) \cap L_m^{Q_1}(\mathbb{R}^2)),$$

$$Z(x,t) \in L_{loc}^1\big([0,\infty)\ ; L_m^{Q_0}(\mathbb{R}^2) \cap L_m^{Q_1}(\mathbb{R}^2)\big),$$

$$\psi(x,t) \in C_0^\infty\big([0,\infty); S_m(\mathbb{R}^2)\big).$$

The definition above corresponds to $m = 1$. It is not hard to see that any weak solution with $m = 2,3,4,\dots$ would also satisfy the condition above (for $m = 1$). This means the integral identity (4.7) holds for all

$$\psi(x,t) \in C_0^\infty\big([0,\infty); S(\mathbb{R}^2)\big).$$

We will term these solutions $R_{\frac{2\pi}{m}}$ -*invariant weak solutions.*

# § 5. Decay of the eigenfunctions of the modified spectral problem

Let $V(x)$ be a radially symmetric velocity field in $\mathbb{R}^2$, such that

$$V(x) = \Omega(|x|)x^{\perp};$$  (5.1)

$$\Omega(|x|) \in C^{\infty}(\mathbb{R}^2);$$  (5.2)

$$\Omega(s) = \kappa^{-1}\left(Cs^{-2} + \frac{1}{2-\alpha}s^{-\alpha}\right), \qquad s \geq M > 0,$$  (5.3)

where $\alpha \in (0,1)$.

An integer $m > 1$ exists so that the following conditions hold true.

Consider the spectral problem

$$L\,\omega \equiv -(V, \partial_y)\omega - (w, \partial_x)\,curl\,V + \alpha\omega + (y, \partial_y)\omega = \lambda\,w\,,$$  (5.4)

where $w = \mathcal{K} * \omega$. Here $\omega \in L^2_m(\mathbb{R}^2)$,

$$\omega(x) = g(|x|)e^{im\theta}\ [i.e., l = 1].$$

Then the spectrum

$$\sigma_{L^2_m(\mathbb{R}^2)}(L) = \{\lambda\} \vee \{\bar{\lambda}\} \vee \left\{\sigma_{L^2_m(\mathbb{R}^2)}(L) \cap \left\{z \in \mathbb{C} | Re\,z < \frac{1}{100}a\right\}\right\}\ (disjoint\ union),$$  (5.5)

where $\lambda$ and $\bar{\lambda}$ are simple eigenvalues of discrete spectrum of $L$ in $L^2_m(\mathbb{R}^2)$, $\lambda = a + ib$; $a > 0, b > 0$.

Notice that $\sigma_{L^2_m(\mathbb{R}^2)}(L) \cap \{z \in \mathbb{C} | Re\,z > \alpha - 1\}$ is a discrete subset of the half-plane $\{z \in \mathbb{C} | Re\,z > \alpha - 1\}$, consisting of eigenvalues with finite multiplicity as in the Riesz-Schauder theory.

To construct such a vector-field $V$, we consider the problem (3.1) with a small parameter $\kappa > 0$, where

$$R(s) = C\,s^{-2} + \frac{1}{2-\alpha}\,s^{-\alpha}, s \geq M > 0.$$

Choose $d = \frac{1}{100}Re\,\lambda_0$ in the Proposition 3.1 and construct $\kappa_0$ as in the Proposition 3.1. Replace $V$ by $\kappa_0^{-1}V$,

$$\lambda = \frac{1}{\kappa_0}\lambda_0 - 1 + \alpha.$$

Since

$$\frac{1}{100}\left(\frac{1}{\kappa_0}Re\,\lambda_0 - 1 + \alpha\right) > \frac{d}{\kappa_0} - 1 + \alpha,$$

(5.5) follows from the Proposition 3.1 for sufficiently small $\kappa_0 > 0$.

We may also assume

$$Re\,\lambda > \alpha\,, \tag{5.6}$$

by choosing $\kappa_0$ small enough.

Proposition 5.1. The following inequality holds

$$|\omega(x)| \le C(1+|x|)^{-m-\alpha-2}, \qquad \forall x \in \mathbb{R}^2. \tag{5.7}$$

Proof. From (1.7) $\quad \omega \in L_m^2(\mathbb{R}^2),\ \omega(x) = g(|x|)e^{im\theta},\ \text{where } m \ge 2 \text{ implies } \omega \in L_m^\infty(\mathbb{R}^2).$

Therefore, (5.4) yields

$$-im\Omega(s)g(s) + (\alpha - \lambda)g(s) + sg'(s) = f(s), \tag{5.8}$$

where

$$f(s) = \frac{i}{2}G'(s)\left(s^{m-1}\int_s^\infty g(\tau)\,\tau^{1-m}\,d\tau + s^{-1-m}\int_0^s g(\tau)\tau^{1+m}\,d\tau\right). \tag{5.9}$$

Therefore,

$$|f(s)| \le C(1+s)^{-1-\alpha},\ \ \forall s \in (0,\infty). \tag{5.10}$$

Thus,

$$g(s) = Cs^{\lambda-\alpha}\exp\left\{-im\int_s^\infty \frac{\Omega(\tau)}{\tau}d\tau\right\} - s^{\lambda-\alpha}\int_s^\infty \tau^{\alpha-\lambda-1}\exp\left\{-im\int_s^\tau \frac{\Omega(\xi)}{\xi}d\xi\right\}f(\tau)d\tau\,. \tag{5.11}$$

The second integral in (5.11) is $O(s^{-1-\alpha})$ as $s \to \infty$. Therefore, $C = 0$ in (5.11) and

$$|g(s)| \le C(1+s)^{-1-\alpha}, \qquad \forall s \in (0,\infty).$$

From (5.9)

$$|f(s)| \le C(1+s)^{-1-\alpha},\ \ \forall s \in (0,\infty).$$

Continuing the bootstrapping we arrive at the inequality

$$|g(s)| \le C(1+s)^{-m-2-\alpha}, \qquad \forall s \in (0,\infty). \tag{5.12}$$

This concludes the proof of the Proposition 5.1. QED

Using (5.9), (5.12) we get

$$f(s) = \frac{i}{2\kappa_0}(-\alpha)\left(\int_0^\infty g(\tau)\tau^{1+m}\,d\tau\right)s^{-m-\alpha-2} + O(s^{-m-2\alpha-2}),$$

as $s \to \infty$.

Therefore, from (5.11)

$$g(s) = -s^{\lambda-\alpha}\int_s^\infty \tau^{\alpha-\lambda-1}\left(1 + O(s^{-\alpha})\right)\left[\frac{i}{2\kappa_0}(-\alpha)\left(\int_0^\infty g(\xi)\xi^{1+m}\,d\xi\right)\tau^{-m-\alpha-2} + O(\tau^{-m-2\alpha-2})\right]d\tau$$

$$= \frac{i}{2\kappa_0}\alpha s^{\lambda-\alpha}\left[\left(\int_0^\infty g(\xi)\xi^{1+m}\,d\xi\right)\int_s^\infty \tau^{-m-\lambda-3} + O(\tau^{-m-\lambda-3-\alpha})d\tau\right] + O(s^{-m-2\alpha-2}) =$$

$$= \frac{i}{2\kappa_0}\alpha\frac{1}{m+\lambda+2}\left(\int_0^\infty g(\xi)\xi^{1+m}\,d\xi\right)s^{-m-\alpha-2} + O(s^{-m-2\alpha-2}), \qquad (5.13)$$

as $s \to \infty$.

Assuming $\int_0^\infty g(\xi)\xi^{1+m}\,d\xi \neq 0$ and normalizing $g$ so that

$$g(s) \sim s^{-m-\alpha-2},$$

as $s \to \infty$, we get

$$\lambda = \frac{i}{2\kappa_0}\alpha\left(\int_0^\infty g(\xi)\xi^{1+m}\,d\xi\right) - m - 2. \qquad (5.14)$$

Lemma 5.1. If $\int_0^\infty g(\xi)\xi^{1+m}\,d\xi = 0$, then necessarily $g \equiv 0$.

Proof. Assuming $\int_0^\infty g(\xi)\xi^{1+m}\,d\xi = 0$ we define

$$h(s) = \sup_{\tau \geq s}|g(\tau)|(1+\tau)^{m+2\alpha+2}, \quad s \in [0,\infty).$$

Evidently, $h(s)$ is bounded, non-increasing on $[0,\infty)$.

We have because of the assumption

$$f(s) = \frac{i}{2}G'(s)\left(s^{m-1}\int_s^\infty g(\tau)\,\tau^{1-m}\,d\tau - s^{-1-m}\int_s^\infty g(\tau)\tau^{1+m}\,d\tau\right). \qquad (5.15)$$

Therefore,

$$|f(s)| \leq C(1+s)^{-1-\alpha}[s^{m-1}\int_s^\infty h(\tau)(1+\tau)^{-m-2\alpha-2}\tau^{1-m}d\tau$$

$$+ s^{-m-1}\int_s^\infty h(\tau)(1+\tau)^{-m-2\alpha-2}\tau^{1+m}d\tau]$$

$$\leq C(1+s)^{-1-\alpha}h(s)(1+s)^{-m-2\alpha-1}$$

$$\leq C(1+s)^{-m-3\alpha-2}h(s) \quad , \forall s \geq 1.$$

From (5.11), where $C = 0$

$$|g(s)| \leq C s^{Re\lambda-\alpha}\int_s^\infty \tau^{\alpha-Re\lambda-1}(1+\tau)^{-m-3\alpha-2}h(\tau)d\tau$$

$$\leq C s^{Re\lambda-\alpha}(1+s)^{-m-3\alpha-2}h(s)\int_s^\infty \tau^{\alpha-Re\lambda-1}d\tau$$

$$\leq C(1+s)^{-m-3\alpha-2}h(s), \forall s \geq 1. \tag{5.16}$$

Thus, for any $\tau \geq 1$

$$|g(\tau)|(1+\tau)^{m+2\alpha+2} \leq C(1+\tau)^{-\alpha}h(\tau).$$

Taking *sup* over $\tau \geq s$, we arrive at the inequality

$$h(s) \leq C(1+s)^{-\alpha}h(s), \quad \forall s \geq 1.$$

This implies $h(s) = 0$ for $s \geq N$, where $N$ is sufficiently large.

From the first inequality in (5.16) we derive

$$h(s) \leq C\int_s^N h(\tau)d\tau, \quad \forall s \in [\eta, N]$$

for any $\eta \in (0, N)$, where the constant $C$ depends on $\eta$ and $N$. By the Gronwall's lemma, $h(s) \equiv 0$

on any interval $[\eta, N]$, and therefore $g(s) \equiv 0$. QED.

# § 6. Constructing a family of solutions

In this section and in what follows, $\Omega(s)$ stands for the angular velocity:

$$V(x) = \Omega(|x|)x^\perp, x \in \mathbb{R}^2.$$

Here we construct a family of solutions to the problem (4.1). To this end, we need to introduce several auxiliary functions.

Let $\chi(|x|) \in C_0^\infty(\mathbb{R}^2)$ be a radial function such that

$$\begin{cases} \tilde{\chi} = 1 \text{ on } [0,1], \ \tilde{\chi} = 0 \text{ on } [2,\infty), \ \ \tilde{\chi} \in C_0^\infty([0,\infty)); \\ \chi(|x|) = \chi_C(|x|) = \tilde{\chi}(C^{-1}|x|), \qquad \forall x \in \mathbb{R}^2 \end{cases} \tag{6.1}$$

for $C = C_4 > 0$.

The dependence of $\chi$ on $C = C_4$ will be usually suppressed in our notation.

We further introduce a function $\mathcal{R}: [0,\infty) \to (0,\infty)$ and a function $\tau: [0,\infty) \to [0,\infty)$.

Let $\omega(x,t); (x,t) \in \mathbb{R}^2 \times [0,\infty)$ be a solution to the system

$$\begin{cases} \partial_t \omega(x,t) = -(v(x,t),\partial_x)\omega(x,t) + Z(x,t); \ x \in \mathbb{R}^2, \ t \in [0,\infty); \\ v(x,t) = (\mathcal{K}(.) * \omega(.,t))(x); \\ \omega(x,0) = \omega_0(x). \end{cases} \tag{6.2}$$

We start with the following *ansatz*:

$$\omega(x,t) = \mathcal{R}(t)^{-\alpha}\sigma\big(y,\tau(t)\big)\big|_{y=\frac{x}{\mathcal{R}(t)}} + \mathcal{R}(t)^{-\alpha}curl_y[\chi(\varepsilon\mathcal{R}(t)|y|)V(y)]\big|_{y=\frac{x}{\mathcal{R}(t)}}. \tag{6.3}$$

Here $\varepsilon > 0$ is a parameter that will be of use later.

Differentiating (6.3) yields

$$\partial_t \omega(x,t) = \dot{\mathcal{R}}(t)\mathcal{R}(t)^{-1-\alpha}\big(-\alpha\sigma(y,\tau(t)) - (y,\partial_y)\sigma(y,\tau(t))\big)\big|_{y=\frac{x}{\mathcal{R}(t)}}$$

$$+\dot{\tau}(t)\mathcal{R}(t)^{-\alpha}\partial_\tau \sigma\big(y,\tau(t)\big)\big|_{y=\frac{x}{\mathcal{R}(t)}} + Z(x,t). \tag{6.4}$$

Here

$$Z(x,t) = \partial_t[\mathcal{R}(t)^{-\alpha} \chi(\varepsilon|x|)curl_y V\left(\frac{x}{\mathcal{R}(t)}\right) + \varepsilon|x|\chi'(\varepsilon|x|)\mathcal{R}(t)^{-\alpha}\Omega\left(\frac{|x|}{\mathcal{R}(t)}\right)]$$

$$= \chi(\varepsilon|x|)\dot{\mathcal{R}}(t)\mathcal{R}(t)^{-1-\alpha}W\left(\frac{x}{\mathcal{R}(t)}\right) + \varepsilon|x|\chi'(\varepsilon|x|)\dot{\mathcal{R}}(t)\mathcal{R}(t)^{-1-\alpha}B\left(\frac{|x|}{\mathcal{R}(t)}\right), \tag{6.5}$$

where

$$\begin{cases} W(y) = -\alpha \ curlV(y) - \big(y, \partial_y\big)curlV(y); \\ B(|y|) = -\alpha \ \Omega(|y|) - |y|\Omega'(|y|). \end{cases}$$ (6.6)

From (1.2) and (6.6),

$$W(y) = 0 \ for \ |y| \geq M > 0,$$

$$B(|y|) = C|y|^{-2} \ for \ |y| \geq M > 0$$

With an appropriate constant $C$.

Therefore, if

$$\mathcal{R}(t) \leq C_4 M^{-1} \varepsilon^{-1},$$ (6.7)

where the constant $C_4$ is coming from (6.1),

$$Z(x,t) = \dot{\mathcal{R}}(t)\mathcal{R}(t)^{-1-\alpha}\left( W\left(\frac{x}{\mathcal{R}(t)}\right) + \varepsilon|x|\chi'(\varepsilon|x|)B\left(\frac{|x|}{\mathcal{R}(t)}\right)\right).$$ (6.8)

We now specify the functions $\mathcal{R}(t), \tau(t)$ and $\sigma(y, \tau)$.

The functions $\mathcal{R}(t)$ and $\tau(t)$ are solutions to the following Cauchy problem:

$$\begin{cases} \dot{\mathcal{R}}(t) = \gamma \ \mathcal{R}(t)^{1-\alpha}; \ \mathcal{R}(0) = 1, \\ \dot{\tau}(t) = \mathcal{R}(t)^{-\alpha}; \ \tau(0) = 0, \end{cases}$$ (6.9)

$\gamma > 0$ is a parameter.

It is not hard to solve the Cauchy problem (6.9):

$$\begin{cases} \mathcal{R}(t) = (1 + \alpha\gamma t)^{\frac{1}{\alpha}}, \quad t \in [0, \infty), \\ \tau(t) = \dfrac{1}{\alpha\gamma}\log(1 + \alpha\gamma t), \quad t \in [0, \infty). \end{cases}$$ (6.10)

Also,

$$\mathcal{R}(t) = e^{\gamma\tau(t)}, t \in [0, \infty).$$ (6.11)

The function $\sigma(y, \tau)$ is defined as the solution of the following Cauchy problem:

$$\begin{cases} \partial_\tau\sigma(y,\tau) = -\big(\chi(\varepsilon\mathcal{R}(t(\tau))|y|)V(y), \partial_y\big)\sigma(y,\tau) - \big(w(y,\tau), \partial_y\big)curl_y[\chi(\varepsilon\mathcal{R}(t(\tau))|y|)V(y)] \\ \qquad\qquad -(w, \partial_y)\sigma + \gamma\big(\alpha\sigma + (y, \partial_y)\sigma\big); \\ \qquad\qquad w(y,\tau) = \big(\mathcal{K}(.) * \sigma(.,\tau)\big)(y); \\ \qquad\qquad \sigma(y,0) = \sigma_0(y). \end{cases}$$ (6.12)

Here $\mathcal{R}(t(\tau)) = e^{\gamma\tau}$, $\mathcal{K}(.)$ denotes the 2-dimensional Biot-Savart kernel.

The following remark is of central importance to our construction.

Theorem 6.1. Given the ansatz (6.3) and (6.5), (6.9), (6.10), (6.12) the vorticity $\omega(x,t)$; $(x,t) \in \mathbb{R}^2 \times [0,\infty)$ satisfies (6.2). The initial condition $\omega_0(x)$ is defined from (6.3):

$$\omega_0(x) = \sigma_0(x) + curl_x[\chi(\varepsilon|x|)V(x)].$$

Proof. We use (6.12), (6.3), (6.4), (6.9) to get

$$\partial_t \omega(x,t) = \gamma\mathcal{R}(t)^{-2\alpha}(-\alpha\sigma(y,\tau(t)) - (y,\partial_y)\sigma(y,\tau(t))|_{y=\frac{x}{\mathcal{R}(t)}}$$

$$+ \mathcal{R}(t)^{-2\alpha}[-(\chi(\varepsilon\mathcal{R}(t(\tau))|y|)V(y),\partial_y)\sigma(y,\tau) - (w(y,\tau),\partial_y)curl_y[\chi(\varepsilon\mathcal{R}(t(\tau))|y|)V(y)]$$

$$-(w,\partial_y)\sigma + \gamma(\alpha\sigma + (y,\partial_y)\sigma)]|_{y=\frac{x}{\mathcal{R}(t)}} + Z(x,t). \qquad (6.13)$$

From (6.13)

$$\partial_t \omega(x,t) = \mathcal{R}(t)^{-2\alpha}[-(\chi(\varepsilon\mathcal{R}(t(\tau))|y|)V(y),\partial_y)\sigma(y,\tau) - (w(y,\tau),\partial_y)curl_y[\chi(\varepsilon\mathcal{R}(t(\tau))|y|)V(y)]$$

$$-(w,\partial_y)\sigma]|_{y=\frac{x}{\mathcal{R}(t)}} + Z(x,t)$$

$$= ((\mathcal{R}(t)^{1-\alpha}\chi(\varepsilon|x|)V\left(\frac{x}{\mathcal{R}(t)}\right) + \mathcal{R}(t)^{1-\alpha}w\left(\frac{x}{\mathcal{R}(t)},\tau(t)\right),\partial_x)$$

$$(\mathcal{R}(t)^{-\alpha}curl_y[\chi(\varepsilon\mathcal{R}(t(\tau))|y|)V(y)]|_{y=\frac{x}{\mathcal{R}(t)}} + \mathcal{R}(t)^{-\alpha}\sigma\left(\frac{x}{\mathcal{R}(t)},\tau(t)\right))$$

$$+ Z(x,t).$$

But

$$curl_x\left[\mathcal{R}(t)^{1-\alpha}\chi(\varepsilon|x|)V\left(\frac{x}{\mathcal{R}(t)}\right)\right] = \mathcal{R}(t)^{-\alpha}curl_y[\chi(\varepsilon\mathcal{R}(t(\tau))|y|)V(y)]|_{y=\frac{x}{\mathcal{R}(t)}}.$$

Therefore, using (6.3) again, we get

$$\partial_t \omega(x,t) = -\left((\mathcal{K}(.) * \omega(.,t))(x),\partial_x\right)\omega(x,t) + Z(x,t)$$

$$= -(v(x,t),\partial_x)\omega(x,t) + Z(x,t).$$

This completes the proof of the Theorem.

It remains to specify $\sigma_0(y)$ in (6.12).

For sufficiently small $\gamma > 0$ we denote by $L$ the following operator

$$L\sigma \equiv -(V(y), \partial_y)\sigma - (w, \partial_y)\,curl\,V(y) + \gamma(\alpha\sigma + (y, \partial_y)\sigma), \tag{6.15}$$

where, as above $w = \mathcal{K} * \sigma$, acting in $L^2_m(\mathbb{R}^2)$, $m \geq 2$ being a fixed integer. The operator $L$ is a closed unbounded operator in $L^2_m(\mathbb{R}^2)$. As was shown above, we can select $V \in \mathcal{C}$, so that the spectrum of the operator $L$ in $L^2_m(\mathbb{R}^2)$ has the following properties:

$$\sigma_{L^2_m(\mathbb{R}^2)}(L) = \{\lambda\} \cup \{\bar{\lambda}\} \cup (\sigma_{L^2_m(\mathbb{R}^2)}(L) \cap \{z \in \mathbb{C} \,|\, Re\,z < \tfrac{a}{100}\}), \tag{6.16}$$

where $\lambda, \bar{\lambda}$ are simple eigenvalues of the discrete spectrum of $L$ in $L^2_m(\mathbb{R}^2)$, i.e., they satisfy the conditions of the Riesz-Schauder theory.

Here,

$$\lambda = a + bi \,; a > 0, b > 0. \tag{6.17}$$

Let $\eta \in L^2_m(\mathbb{R}^2)$ be a normalized solution to the spectral problem

$$\begin{cases} L\,\eta = \lambda\,\eta; \\ \eta(x) = g(|x|)e^{im\theta}; \\ g(s) \sim s^{-m-\alpha-2}, \quad as\ s \to \infty. \end{cases} \tag{6.18}$$

We define

$$\sigma_0(y) = Re\left(\varepsilon^{\rho+i\zeta}\,\eta(y)\right), \quad y \in \mathbb{R}^2; \quad \rho, \zeta \in \mathbb{R}\,, \quad \rho > 0\,. \tag{6.19}$$

From the Proposition 5.1 (see (5.7))

$$\sigma_0 \in L^1_m(\mathbb{R}^2) \cap L^\infty_m(\mathbb{R}^2). \tag{6.20}$$

Remark. For the eigenfunction $\eta$ such as in (6.18) we have

$$|y|\partial_y\eta(y) \in L^1_m(\mathbb{R}^2) \cap L^\infty_m(\mathbb{R}^2). \tag{6.21}$$

Indeed, this follows from (6.20) and (5.8) - (5.9).

Remark. We will need

$$\gamma < \frac{a}{\alpha} \tag{6.22}$$

below.

# § 7. Controlling the time evolution as $\tau \to \infty$.

In this section we study the growth of solutions to the Euler-type equations as $\tau \to \infty$.

To simplify the notation in this section, we often suppress the index $m \in \mathbb{Z}$, $m \geq 2$ of a function space.

Let $\sigma(y, \tau)$ be a solution to the following problem:

$$\begin{cases} \partial_\tau \sigma(y, \tau) = -(V(y), \partial_y)\sigma(y, \tau) - (w(y, \tau), \partial_y)curl_y[V(y)] \\ \qquad -(w, \partial_y)\sigma + \gamma(\alpha\sigma + (y, \partial_y)\sigma); \\ \qquad w(y, \tau) = (\mathcal{K}(.) * \sigma(., \tau))(y); \\ \qquad \sigma(y, 0) = Re\ (\varepsilon^{\rho + i\zeta}\eta(y)), \qquad y \in \mathbb{R}^2. \end{cases} \tag{7.1}$$

Here $\gamma > 0$ is a sufficiently small constant, so that the partition (6.16) of the spectrum of $L_m^2$ holds;

$$\rho, \zeta \in \mathbb{R}, \quad \rho > 0$$

will be selected later.

We choose $Q \in \left(2, \frac{2}{\alpha}\right)$.

Let $C_2 > 0$ be a sufficiently large constant that depends only on $V(y)$ so that

$$\begin{cases} \|\eta\|_{L^Q(\mathbb{R}^2)} \leq C_2; \\ \||y|\partial_y\eta\|_{L^Q(\mathbb{R}^2)} \leq C_2; \\ \|\partial_y\eta\|_{L^{2,1}(\mathbb{R}^2)} \leq C_2. \end{cases} \tag{7.2}$$

We select a constant $C_1 > CC_2$ to be specified later. Here $C >> 1$ is a sufficiently large constant.

Let $T \in (0, \infty]$ be the *maximal time* so that the solution $\sigma(y, \tau)$ of the equation (7.1), with

$$\sigma(y, 0) = Re\ (\varepsilon^{\rho + i\zeta}\eta(y)), \qquad y \in \mathbb{R}^2 \tag{7.3}$$

satisfies the following three inequalities:

$$\begin{cases} \|\sigma(\tau)\|_{L^Q(\mathbb{R}^2)} \leq C_1\varepsilon^\rho e^{a\tau}, \qquad \tau \in [0, T]; \\ \||y|\partial_y\sigma(\tau)\|_{L^Q(\mathbb{R}^2)} \leq C_1\varepsilon^\rho e^{a\tau}, \qquad \tau \in [0, T]; \\ \|\partial_y\sigma(\tau)\|_{L^{2,1}(\mathbb{R}^2)} \leq C_1\varepsilon^\rho e^{a\tau}, \qquad \tau \in [0, T]. \end{cases} \tag{7.4}$$

We define $T_*$ so that

$$C_1 \varepsilon^\rho e^{aT_*} = C_0 \ll 1, \tag{7.5}$$

where $C_0$ is a fixed positive constant to be chosen later.

Here is the main result of this section. It turns out, $T$ is rather "large".

Theorem 7.1. The following inequality holds:

$$T \geq T_* \quad . \tag{7.6}$$

Proof. We assume to the contrary that $T < T_*$.

We need to estimate $\|\sigma(\tau)\|_{L^Q(\mathbb{R}^2)}$, $\||y|\partial_y \eta\|_{L^Q(\mathbb{R}^2)}$, $\|\partial_y \eta\|_{L^{2,1}(\mathbb{R}^2)}$. The spatial derivatives of $\sigma$ in (7.1) satisfy the following system of equations:

$$\partial_\tau \partial_k \sigma(y, \tau) = -(V(y), \partial_y)\partial_k \sigma(y, \tau) - (\partial_k V(y), \partial_y)\sigma$$

$$+\gamma\big((\alpha + 1)\partial_k \sigma + (y, \partial_y)\partial_k \sigma\big)$$

$$-(\partial_k w(y, \tau), \partial_y)curl_y[V(y)] - (w(y, \tau), \partial_y)\partial_k curl_y[V(y)]$$

$$-(w, \partial_y)\partial_k \sigma - (\partial_k w, \partial_y)\sigma; \qquad \partial_k \sigma(y, 0) = \partial_k \sigma_0(y), \ \ k = 1, 2. \tag{7.7}$$

Here as above

$$w(y, \sigma) = \big(\mathcal{K}(.) * \sigma(., \tau)\big)(y).$$

From the variation of parameters formula, applied to (7.1), for any $\tau \in [0, \infty)$

$$\sigma(\tau) = e^{\tau L}\sigma_0 + \int_0^\tau e^{(\tau - \xi)L}\left[-\big(w(\xi), \partial_y\big)\sigma(\xi)\right]d\xi$$

in $L^Q(\mathbb{R}^2)$. Therefore,

$$\|\sigma(\tau)\|_{L^Q(\mathbb{R}^2)} \ \leq \ e^{a\tau}\|\sigma_0\|_{L^Q(\mathbb{R}^2)} + C\int_0^\tau e^{a(\tau - \xi)}\|\big(w(\xi), \partial_y\big)\sigma(\xi)\|_{L^Q(\mathbb{R}^2)}d\xi \ . \tag{7.8}$$

But from (7.4) for any $\xi \in [0, \tau] \subset [0, T]$

$$\|\big(w(\xi), \partial_y\big)\sigma(\xi)\|_{L^Q(\mathbb{R}^2)} \leq \ \|\frac{w(\xi)}{|y|}\|_{L^\infty(\mathbb{R}^2)} \quad \||y|\partial_y \sigma(\xi)\|_{L^Q(\mathbb{R}^2)}$$

$$\leq C\|\partial_y w(\xi)\|_{L^\infty(\mathbb{R}^2)}\||y|\partial_y \sigma(\xi)\|_{L^Q(\mathbb{R}^2)}$$

$$\leq C\|\partial_y \sigma(\xi)\|_{L^{2,1}(\mathbb{R}^2)}\||y|\partial_y \sigma(\xi)\|_{L^Q(\mathbb{R}^2)}$$

$$\leq CC_1^2 \varepsilon^{2\rho} e^{2a\xi}. \tag{7.9}$$

Therefore, using (7.8), (7.9) for any $\tau \in [0, T]$

$$\|\sigma(\tau)\|_{L^Q(\mathbb{R}^2)} \leq C_2 \varepsilon^\rho e^{\alpha\tau} + CC_1^2 \varepsilon^{2\rho} \int_0^\tau e^{\alpha\tau + a\xi}\, d\xi$$

$$\leq C_2 \varepsilon^\rho e^{\alpha\tau} + CC_1^2 \varepsilon^{2\rho} e^{2a\tau}$$

$$\leq C_2 \varepsilon^\rho e^{\alpha\tau} + CC_0 C_1 \varepsilon^\rho e^{\alpha\tau}, \tag{7.10}$$

where on the last step we used the assumption $T < T_*$ (and $\tau \in [0, T]$).

If $C_0$ is chosen so that in the right side of (7.10)

$$CC_0 C_1 \leq C_2, \tag{7.11}$$

Then

$$\|\sigma(\tau)\|_{L^Q(\mathbb{R}^2)} \leq 2C_2 \varepsilon^\rho e^{\alpha\tau}, \qquad \forall\, \tau \in [0, T]. \tag{7.12}$$

Next, we estimate $\|\partial_y \sigma(\tau)\|_{L^{2,1}(\mathbb{R}^2)}$, using (7.7).

Lemma 7.1.

Let

$$\begin{cases} \dfrac{dx}{d\tau} = V(x) + w(x, \tau) - \gamma x; \\ \quad x(0) = y \in \mathbb{R}^2. \end{cases}$$

Then, for any $\tau \in [0, T]$ and given that $C_0$ is sufficiently small, the following inequality holds:

$$\int_0^\tau |\frac{1}{2}\left(\left(\partial_y V(x(\xi))\right) + \left(\partial_y V(x(\xi))\right)^*\right)|d\xi \leq C < \infty, \tag{7.13}$$

Where the constant $C$ in (7.13) is uniform with respect to $y \in \mathbb{R}^2$. Here $(\partial_y V)^*$ stands for the transposition of the matrix $\partial_y V$.

Proof. We have

$$\frac{d}{d\tau}|x(\tau)| = \left(\frac{w(x(\tau),\tau)}{|x(\tau)|}, x(\tau)\right) - \gamma|x(\tau)|;$$ (7.14)

$$\|\frac{w(.,\tau)}{|.|}\|_{L^\infty(\mathbb{R}^2)} \le C\|\partial_y w(\tau)\|_{L^\infty(\mathbb{R}^2)}$$

$$\le C\|\partial_y \sigma(\tau)\|_{L^{1,2}(\mathbb{R}^2)} \le CC_1 \varepsilon^\rho e^{\alpha\tau}, \ \ \forall \tau \in [0,T].$$ (7.15)

From (7.14), (7.15) we derive the following two inequalities

$$|y|\exp\{-\gamma\tau - CC_1\varepsilon^\rho e^{\alpha\tau}\} \le |x(\tau)| \le |y|\exp\{-\gamma\tau + CC_1\varepsilon^\rho e^{\alpha\tau}\}.$$ (7.16)

Also, $V(y) = \Omega(|y|)y^\perp$, so

$$\frac{1}{2}\left(\left(\partial_y V(y)\right) + \left(\partial_y V(y)\right)^*\right) = \frac{1}{|y|}\ \Omega'(|y|)\begin{pmatrix} -y_1 y_2 & \frac{1}{2}(y_1^2 - y_2^2) \\ \frac{1}{2}(y_1^2 - y_2^2) & y_1 y_2 \end{pmatrix},$$

therefore,

$$\left|\frac{1}{2}\left(\left(\partial_y V(y)\right) + \left(\partial_y V(y)\right)^*\right)\right| \le C|y|^2(1+|y|)^{-2-\alpha}.$$

If $C_0$ is chosen to be small enough, so that $CC_1\varepsilon^\rho e^{\alpha\tau} \le CC_0 \le 1$ in (7.16), then the integral in the left side of (7.13) is

$$\le C\int_0^\tau |ye^{-\gamma\xi}|^2 \ \left(1 + e^{-\gamma\xi}|y|\right)^{-2-\alpha} d\xi$$

$$= \frac{1}{\gamma}C\int_{e^{-\gamma\tau}|y|}^{|y|} r \ (1+r)^{-2-\alpha} dr$$

$$\le \frac{1}{\gamma}C\int_0^\infty r \ (1+r)^{-2-\alpha} dr.$$

Since the last integral converges, this completes the proof of the lemma. QED.

The following estimate and its analogue below for Lorentz spaces will be used frequently in the sequel. We present a general scheme, in applying the lemma we will deal with classical $C^\infty$ solutions where all the computations in the proof can be easily justified.

Lemma 7.2. Let $\mu, \nu \in \mathbb{R}$. Let $W(y,\tau) = \left(W_1(y,\tau), ..., W_N(y,\tau)\right)$ be an $N-$vector, $y \in \mathbb{R}^n, \tau \in [0,T]$,

satisfying the following system of equations

$$\begin{cases} \partial_\tau W(y,\tau) = -(Z(y,\tau),\partial_y)W - A(y,\tau)W - \mu W + F(y,\tau) \\ W(y,0) = W_0(y) \in (L^q(\mathbb{R}^n))^N, q \in [1,\infty]. \end{cases} \quad (7.17)$$

Assume $Z(y,\tau) = (Z_1(y,\tau), \dots, Z_n(y,\tau))$ be continuous in $(y,\tau)$ and Lipschitz continuous in $y$ variable

$$\|\partial_y Z(\tau)\|_{L^\infty(\mathbb{R}^n)} \le C < \infty, \ \ \forall \, \tau \in [0,T];$$

$$A(y,\tau) \in Mat_{n \times n}(\mathbb{R});$$

$$\|A(\tau)\|_{L^\infty(\mathbb{R}^n)} \le C < \infty, \ \ \forall \, \tau \in [0,T]$$

Let $x(\tau)$ be the unique flow, generated by $Z(y,\tau)$:

$$\frac{d}{d\tau} x(\tau) = Z(x(\tau),\tau); \ \ x(0) = y \in \mathbb{R}^n.$$

Assume the following estimates:

$$\|div_y Z(\tau) - \nu\|_{L^\infty(\mathbb{R}^n)} \le g(\tau), \ \ \forall \, \tau \in [0,T];$$

$$\int_\zeta^\tau \frac{1}{2} |(A + A^*)(x(\xi),\xi)| d\xi \le H(\tau,\zeta) \ \ \ \ \forall \, \zeta, \tau \in [0,T], \ \ \zeta \le \tau.$$

Assume $F \in L^1([0,T]; (L^q(\mathbb{R}^n))^N)$. Then

$$\|W(\tau)\|_{L^q(\mathbb{R}^n)} \le \|W(0)\|_{L^q(\mathbb{R}^n)} \exp\{\int_0^\tau \frac{1}{q} g(\xi) d\xi + H(\tau,0) + (\frac{\nu}{q} - \mu)\tau\}$$

$$+ \int_0^\tau \|F(\xi)\|_{L^q(\mathbb{R}^n)} \exp\{\int_\xi^\tau \frac{1}{q} g(\zeta) d\zeta + H(\tau,\xi) + (\frac{\nu}{q} - \mu)(\tau - \xi)\} \, d\xi.$$

Proof. By Duhamel's principle it is sufficient to prove Lemma with $F = 0$.

Differentiating along the trajectory $x(\tau)$ of the flow generated by the vector field $Z(x,\tau)$ and using (7.17) yields

$$\frac{d}{d\tau} W(x(\tau),\tau) = -A(x(\tau),\tau)W(x(\tau),\tau) - \mu \, W(x(\tau),\tau).$$

Therefore,

$$\frac{d}{d\tau}|W(x(\tau),\tau)|^2 = -\big((A+A^*)(x(\tau),\tau)W(x(\tau),\tau),W(x(\tau),\tau)\big) - 2\mu\,|W(x(\tau),\tau)|^2;$$

$$|W(x(\tau),\tau)|^2 \le \exp\{2H(\tau,0) - 2\,\mu\,\tau\}\,|W(y,0)|^2.$$

Changing variables from $x(\tau)$ to $y$ yields

$$|W(x(\tau),\tau)|^q|dx(\tau)| = |W(x(\tau),\tau)|^q\det\left(\frac{\partial x(\tau)}{\partial y}\right)|dy|$$

$$\le |W(y,0)|^q\exp\{qH(\tau,0) - q\mu\tau\}\,\exp\{\nu\tau + \int_0^\tau g(\xi)d\xi\}\,|dy|,$$

using Liouville's theorem. Integrating in $y$, we arrive at the estimate

$$\|W(\tau)\|_{L^q(\mathbb{R}^n)} \le \|W(0)\|_{L^q(\mathbb{R}^n)}\exp\left\{\int_0^\tau\frac{1}{q}g(\xi)d\xi + H(\tau,0) + \left(\frac{\nu}{q} - \mu\right)\tau\right\}.$$

This completes the proof of the Lemma. QED.

We need similar estimates in Lorentz spaces, given in the following Lemma.

Lemma 7.3. Let $\mu,\nu\in\mathbb{R}$. Let $W(y,\tau) = \big(W_1(y,\tau),\dots,W_N(y,\tau)\big)$ be an $N$-vector, $y\in\mathbb{R}^n,\tau\in[0,T]$,

satisfying the following system of equations

$$\begin{cases}\partial_\tau W(y,\tau) = -\big(Z(y,\tau),\partial_y\big)W - A(y,\tau)W - \mu W + F(y,\tau)\\ \qquad W(y,0) = W_0(y)\in(L^{q,1}(\mathbb{R}^n))^N, q\in(1,\infty)\end{cases}$$

Assume $Z(y,\tau) = \big(Z_1(y,\tau),\dots,Z_n(y,\tau)\big)$ be continuous in $(y,\tau)$ and Lipschitz continuous in $y$ variable

$$\|\partial_y Z(\tau)\|_{L^\infty(\mathbb{R}^n)} \le C < \infty,\ \ \forall\,\tau\in[0,T];$$

$$A(y,\tau)\in Mat_{n\times n}(\mathbb{R});$$

$$\|A(\tau)\|_{L^\infty(\mathbb{R}^n)} \le C < \infty,\ \ \forall\,\tau\in[0,T]$$

Let $x(\tau)$ be the unique flow, generated by $Z(y,\tau)$:

$$\frac{d}{d\tau}x(\tau) = Z(x(\tau),\tau);\ \ x(0) = y\in\mathbb{R}^n.$$

Assume the following estimates:

$$\|div_y Z(\tau) - \nu\|_{L^\infty(\mathbb{R}^n)} \leq g(\tau), \quad \forall \, \tau \in [0, T];$$

$$\int_\zeta^\tau \frac{1}{2} |(A + A^*)(x(\xi), \xi)| d\xi \leq H(\tau, \zeta) \quad \forall \, \zeta, \tau \in [0, T], \quad \zeta \leq \tau.$$

Assume $F \in L^1([0, T]; (L^{q,1}(\mathbb{R}^n))^N)$.

Then

$$\|W(\tau)\|_{L^{q,1}(\mathbb{R}^n)} \leq C \|W(0)\|_{L^{q,1}(\mathbb{R}^n)} \exp\{\int_0^\tau \frac{1}{q} g(\xi) d\xi + H(\tau, 0) + (\frac{\nu}{q} - \mu)\tau\}$$

$$+ C \int_0^\tau \|F(\xi)\|_{L^{q,1}(\mathbb{R}^n)} \exp\{\int_\xi^\tau \frac{1}{q} g(\zeta) d\zeta + H(\tau, \xi) + (\frac{\nu}{q} - \mu)(\tau - \xi)\} \, d\xi.$$

Here a constant $C$ depends only on $n, N, q$ but not on the initial data $W_0 \in (L^{q,1}(\mathbb{R}^n))^N$.

We will be using the following norm on $L^{q,1}(\mathbb{R}^n)$, $q \in (1, \infty)$. Let for a measurable function $f$ on $\mathbb{R}^n$

$$f^{**}(t) = \sup_{|E| \geq t} \frac{1}{|E|} \int_E |f| dx, \quad \forall \, t \in (0, \infty);$$

The supremum is taken over all Borel measurable sets $E \subset \mathbb{R}^n$ of finite measure $|E| \geq t$. We define

$$\|f\|_{L^{q,1}(\mathbb{R}^n)} = \int_0^\infty t^{\frac{1}{q}} f^{**}(t) \frac{dt}{t}.$$

Let also $\forall j \in \mathbb{Z}$ $E_j = \{x \in \mathbb{R}^n | 2^j \leq |f(x)| < 2^{j+1}\}$. We define

$$\|\|f\|\|_{L^{q,1}(\mathbb{R}^n)} = \sum_{j \in \mathbb{Z}} 2^j |E_j|^{\frac{1}{q}}.$$

Remark. We do not claim $\|\|.\|\|_{L^{q,1}(\mathbb{R}^n)}$ is a norm on $L^{q,1}(\mathbb{R}^n)$.

Proposition 7.1. For any $q \in (1, \infty)$ there exists a constant $C_q > 0$ that depends only on $q$ such that for any measurable $f$

$$C_q^{-1} \|\|f\|\|_{L^{q,1}(\mathbb{R}^n)} \leq \|f\|_{L^{q,1}(\mathbb{R}^n)} \leq C_q \|\|f\|\|_{L^{q,1}(\mathbb{R}^n)}.$$

Proof of Proposition 7.1. First,

$$\|f\|_{L^{q,1}(\mathbb{R}^n)} \leq \sum_{j \in \mathbb{Z}} \|f 1_{E_j}\|_{L^{q,1}(\mathbb{R}^n)}$$

$$\leq \sum_{j \in \mathbb{Z}} 2^{j+1} \|1_{E_j}\|_{L^{q,1}(\mathbb{R}^n)}.$$

Let $g = 1_{E_j}$. Then,

$$g^{**}(t) = \begin{cases} 1, & t \leq |E_j| \\ t^{-1}|E_j|, & t \geq |E_j|. \end{cases}$$

Therefore,

$$\|g\|_{L^{q,1}(\mathbb{R}^n)} = \int_0^{|E_j|} t^{\frac{1}{q}} \frac{dt}{t} + |E_j| \int_{|E_j|}^\infty t^{\frac{1}{q}-1} \frac{dt}{t} = \frac{q^2}{q-1} |E_j|^{\frac{1}{q}}.$$

Thus,

$$\|f\|_{L^{q,1}(\mathbb{R}^n)} \leq 2 \frac{q^2}{q-1} \|f\|_{L^{q,1}(\mathbb{R}^n)}.$$

To prove the opposite inequality, we take $E^{(k)} = \underset{j \geq k}{\cup} E_j = \{x \in \mathbb{R}^n \mid |f(x)| \geq 2^k\}$.

$$\frac{1}{|E^{(k)}|} \int_{E^{(k)}} |f| dx \geq \frac{1}{|E^{(k)}|} \sum_{j \geq k} 2^j |E_j|.$$

Then for $t \in [|E^{(k+1)}|, |E^{(k)}|)$,

$$f^{**}(t) \geq \frac{1}{|E^{(k)}|} \sum_{j \geq k} 2^j |E_j|.$$

Thus,

$$\|f\|_{L^{q,1}(\mathbb{R}^n)} = \int_0^\infty t^{\frac{1}{q}} f^{**}(t) \frac{dt}{t}$$

$$= \sum_{k \in \mathbb{Z}} \int_{|E^{(k+1)}|}^{|E^{(k)}|} t^{\frac{1}{q}} f^{**}(t) \frac{dt}{t} \geq$$

$$\geq q \sum_{k \in \mathbb{Z}} \frac{\sum_{j \geq k} 2^j |E_j|}{\sum_{j \geq k} |E_j|} (|E^{(k)}|^{\frac{1}{q}} - (|E^{(k+1)}|)^{\frac{1}{q}})$$

$$\geq q \sum_{k \in \mathbb{Z}} 2^k \left( |E^{(k)}|^{\frac{1}{q}} - \left( |E^{(k+1)}| \right)^{\frac{1}{q}} \right)$$

$$= \frac{q}{2} \sum_{k \in \mathbb{Z}} 2^k |E^{(k)}|^{\frac{1}{q}} \geq \frac{q}{2} \sum_{k \in \mathbb{Z}} 2^k |E_k|^{\frac{1}{q}} = \| f \|_{L^{q,1}(\mathbb{R}^n)} .$$

To justify the Abel summation formula used above, we need the following statement: $2^k |E^{(k)}|^{\frac{1}{q}} \to 0$ as $k \to \pm\infty$.

Since $|f| \geq 2^k 1_{E^{(k)}}$ and by convergence of the Lebesgue integral $\int_0^\infty t^{\frac{1}{q}} f^{**}(t) \frac{dt}{t}$,

the sequence $\{ 2^k |E^{(k)}|^{\frac{1}{q}} \}_{k \in \mathbb{Z}}$ is bounded. Also, by absolute continuity, as $k \to \infty$

$$2^k \int_0^{|E^{(k)}|} t^{\frac{1}{q}} \frac{dt}{t} \to 0.$$

For $k \to -\infty$ in case we can extract a subsequence $k_m$ so that $|E^{(k_m)}| \to \infty$ as $m \to \infty$,

we get from the absolute continuity of the Lebesgue integral that

$$2^{k_m} |E^{(k_m)}| \int_{|E^{(k_m)}|}^\infty t^{\frac{1}{q}-1} \frac{dt}{t} \to 0,$$

as $m \to \infty$. From here the justification of the Abel summation formula follows. This concludes the Proof of Proposition 7.1. QED.

Proof of Lemma 7.3. By Duhamel's principle it is sufficient to prove Lemma with $F = 0$.

Let $E_j = \{ x \in \mathbb{R}^n | 2^j \leq |W_0(x)| < 2^{j+1} \}$.

We break $W_0(x)$ as follows.

$$W_0(x) = \sum_{j \in \mathbb{Z}} W_0(x) 1_{E_j} \qquad .$$

Differentiating along the trajectory $x(\tau)$ of the flow generated by the vector field $Z(x, \tau)$

$$\frac{d}{d\tau} W(x(\tau), \tau) = -A(x(\tau), \tau) W(x(\tau), \tau) - \mu W(x(\tau), \tau).$$

Therefore,

$$\frac{d}{d\tau} |W(x(\tau), \tau)|^2 = -\big((A + A^*)(x(\tau), \tau) W(x(\tau), \tau), W(x(\tau), \tau)\big) - 2\mu |W(x(\tau), \tau)|^2;$$

$$|W(x(\tau), \tau)|^2 \leq \exp\{2H(\tau, 0) - 2\mu\tau\} |W(y, 0)|^2.$$

Therefore,

$$|W(x(\tau), \tau)| \leq \exp\{H(\tau, 0) - \mu\tau\} |W(y, 0)|$$

Changing variables from $x(\tau)$ to $y$ yields

$$|dx(t)| = \det\left(\frac{\partial x(\tau)}{\partial y}\right) |dy|$$

$$= \exp\left\{\nu\tau + \int_0^\tau g(\xi)d\xi\right\} |dy|,$$

using Liouville's theorem. For $y \in E_j$

$$|W(x(\tau), \tau)| \leq \exp\{H(\tau, 0) - \mu\tau\} 2^{j+1}$$

Therefore,

$$|W(x, \tau)| \leq \exp\{H(\tau, 0) - \mu\tau\} \sum_{j \in \mathbb{Z}} 2^{j+1} 1_{\{x(\tau) | y \in E_j\}} \quad .$$

From this inequality, using Proposition 7.1, we derive

$$\|W(\tau)\|_{L^{q,1}(\mathbb{R}^n)} \leq C \exp\{H(\tau, 0) - \mu\tau\} \exp\left\{q^{-1}\nu\tau + q^{-1}\int_0^\tau g(\xi)d\xi\right\}$$

$$\sum_{j \in \mathbb{Z}} 2^j |E_j|^{\frac{1}{q}}$$

$$\leq C \exp\{H(\tau, 0) - \mu\tau\} \exp\left\{q^{-1}\nu\tau + q^{-1}\int_0^\tau g(\xi)d\xi d\xi\right\} \|W_0\|_{L^{q,1}(\mathbb{R}^n)}.$$

This completes the proof of Lemma 7.3. QED.

Next step is to estimate $\|\partial_y \sigma(\tau)\|_{L^{2,1}(\mathbb{R}^2)}$ using the equation (7.7). We apply Lemma 7.3 to

$$W_k(y,\tau) = \partial_k \sigma(y,\tau); \quad k = 1,2 \, ; N = 2, n = 2;$$

$$Z(y,\tau) = V(y) + w(y,\tau) - \gamma y;$$

$$\mu = -\gamma(\alpha + 1), \nu = -2\gamma;$$

$$A_{kl}(y,\tau) = \partial_k V_l(y) + \partial_k w_l(y,\tau); \quad k,l = 1,2;$$

$$F_k(y,\tau) = -(\partial_k w, \partial_y) \, curl \, V(y) - (w, \partial_y) \partial_k \, curl \, V(y); \quad k = 1,2;$$

$$g(\tau) = 0 \, , \qquad \forall \tau \in [0,T].$$

From Lemma 7.1 and (7.4)

$$\int_0^\tau \frac{1}{2} |(A + A^*)(x(\xi), \xi)| d\xi \le C < \infty \quad , \forall \tau \in [0,T]. \tag{7.18}$$

Indeed,

$$\|\partial_y w(\tau)\|_{L^\infty(\mathbb{R}^2)} \le C \, \|\partial_y \sigma(\tau)\|_{L^{2,1}(\mathbb{R}^2)} \le CC_1 \varepsilon^\rho e^{a\tau};$$

$$\int_0^\tau \|\partial_y w(\xi)\|_{L^\infty(\mathbb{R}^2)} d\xi \le CC_1 \varepsilon^\rho e^{a\tau} \le CC_0 \qquad \forall \tau \in [0,T].$$

To estimate $\|F(\tau)\|_{L^{2,1}(\mathbb{R}^2)}$ we use the following inequalities

$$\|\partial_k w\|_{L^Q(\mathbb{R}^2)} \le C \|\sigma\|_{L^Q(\mathbb{R}^2)} \, , \quad k = 1,2 \, ;$$

$$\|\frac{1}{|y|} w\|_{L^Q(\mathbb{R}^2)} \le C \|\sigma\|_{L^Q(\mathbb{R}^2)}, \quad Q \in (2,\infty) \quad ;$$

$$\begin{cases} \|(\partial_k w, \partial_y) \, curl \, V(y)\|_{L^{2,1}(\mathbb{R}^2)} \le C \|\partial_k w\|_{L^Q(\mathbb{R}^2)} \|\partial_y \, curl \, V\|_{L^{\frac{2Q}{Q-2}, \frac{Q}{Q-1}}(\mathbb{R}^2)} \le C \|\sigma\|_{L^Q(\mathbb{R}^2)} \\ \|(w, \partial_y) \partial_k \, curl \, V\|_{L^{2,1}(\mathbb{R}^2)} \le C \|\frac{1}{|y|} w\|_{L^Q(\mathbb{R}^2)} \||y| \partial_y \partial_k \, curl \, V\|_{L^{\frac{2Q}{Q-2}, \frac{Q}{Q-1}}(\mathbb{R}^2)} \le C \|\sigma\|_{L^Q(\mathbb{R}^2)} \end{cases} \tag{7.19}$$

From (7.19), (7.12)

$$\|F(\tau)\|_{L^{2,1}(\mathbb{R}^2)} \le CC_2 \varepsilon^\rho e^{a\tau}, \quad \forall \tau \in [0,T]. \tag{7.20}$$

Using Lemma 7.3 and the inequalities (7.18) --(7.20) we get

$$\|\partial_y \sigma(\tau)\|_{L^{2,1}(\mathbb{R}^2)} \leq CC_2 \varepsilon^\rho e^{\alpha\gamma\tau} + CC_2 \int_0^\tau e^{\alpha\gamma(\tau-\xi)} \varepsilon^\rho e^{a\xi}\, d\xi \leq CC_2 \varepsilon^\rho e^{a\tau}, \quad \forall \tau \in [0,T], \qquad (7.21)$$

provided $\gamma < \frac{1}{\alpha} a$.

We now turn to the estimate of $\||y|\partial_y \sigma(\tau)\|_{L^Q(\mathbb{R}^2)}$.

From (7.7) we derive the following equation (summation with respect to the indexes repeated twice is assumed):

$$\partial_\tau(|y|\partial_k\sigma(y,\tau)) = -(V, \partial_y)(|y|\partial_k\sigma) - \partial_k V_j\, |y|\partial_j\sigma$$

$$+\gamma(\alpha|y|\partial_k\sigma + (y, \partial_y)(|y|\partial_k\sigma))$$

$$-\partial_k w_j |y|\partial_j curl\, V - w_j |y|\partial_j\partial_k curl\, V$$

$$-(w, \partial_y)(|y|\partial_k\sigma) + |y|^{-2}(w, y)(|y|\partial_k\sigma)$$

$$-\partial_k w_j(|y|\partial_j\sigma). \qquad (7.22)$$

We use Lemma 7.2 with $q = Q$,

$$W_k(y,\tau) = |y|\, \partial_k\sigma(y,\tau); \ \ k = 1,2\, ; N = 2, n = 2;$$

$$Z(y,\tau) = V(y) + w(y,\tau) - \gamma y;$$

$$\mu = -\alpha\gamma, \nu = -2\gamma;$$

$$A_{kj}(y,\tau) = \partial_k V_j(y) + \partial_k w_j(y,\tau) - |y|^{-2}(w,y)\delta_{kj}; \ \ k,j = 1,2;$$

$$g(\tau) = 0, \qquad \forall \tau \in [0,T];$$

$$F_k(y,\tau) = -\partial_k w_j |y|\partial_j curl\, V - w_j |y|\partial_j\partial_k curl\, V, \ \ k = 1,2.$$

We first estimate $\|F(\tau)\|_{L^Q(\mathbb{R}^2)}$. We have

$$\|F(\tau)\|_{L^Q(\mathbb{R}^2)} \leq \|\partial_y w(\tau)\|_{L^Q(\mathbb{R}^2)}\ \||y|\partial_y\, curl\, V\|_{L^\infty(\mathbb{R}^2)}$$

$$+\||y|^{-1} w(\tau)\|_{L^Q(\mathbb{R}^2)}\||y|^2\partial_y^2\, curl\, V\|_{L^\infty(\mathbb{R}^2)}$$

$$\leq C\|\sigma(\tau)\|_{L^Q(\mathbb{R}^2)}. \qquad (7.23)$$

The integral

$$\int_0^\tau \frac{1}{2}|(A + A^*)(x(\xi),\xi)|d\xi$$

contains one additional term as compared to the previous inequality (7.18). This term is bounded above by

$$\int_0^\tau \||y|^{-1} w(y, \xi)\|_{L^\infty(\mathbb{R}^2)}\, d\xi$$

$$\leq C \int_0^\tau \|\partial_y w(\xi)\|_{L^\infty(\mathbb{R}^2)}\, d\xi$$

$$\leq C \int_0^\tau \|\partial_y \sigma(\xi)\|_{L^{2,1}(\mathbb{R}^2)}\, d\xi$$

$$\leq C \int_0^\tau C C_2 \varepsilon^\rho e^{a\xi}\, d\xi$$

$$\leq C C_2 \varepsilon^\rho e^{a\tau} \leq C C_0 C_2 C_1^{-1} \leq C C_0, \quad \forall \tau \in [0, T].$$

Therefore,

$$\int_0^\tau \frac{1}{2}|(A + A^*)(x(\xi), \xi)|\, d\xi \leq C < \infty, \qquad \forall \tau \in [0, T]. \tag{7.24}$$

From (7.22) - (7.24) , (7.12) and Lemma 7.2 we get the following inequality:

$$\||y|\partial_y \sigma(\tau)\|_{L^Q(\mathbb{R}^2)} \leq$$

$$\leq C C_2 \varepsilon^\rho e^{\gamma\left(\alpha - \frac{2}{Q}\right)\tau} + C \int_0^\tau e^{\gamma\left(\alpha - \frac{2}{Q}\right)(\tau - \xi)} \|\sigma(\xi)\|_{L^Q(\mathbb{R}^2)} d\xi$$

$$\leq C C_2 \varepsilon^\rho e^{\gamma\left(\alpha - \frac{2}{Q}\right)\tau} + C C_2 \varepsilon^\rho e^{a\tau}$$

$$\leq C C_2 \varepsilon^\rho e^{a\tau}, \quad \forall \tau \in [0, T]\,. \tag{7.25}$$

We go back now to the definition of $T$ in (7.4). For at least one of the three inequalities in (7.4) there should be an equality as $\tau = T$.

From (7.12), (7.21), (7.25) we derive

$$C_1 \leq C C_2.$$

But this inequality is not possible if $C_1$ has initially been chosen to satisfy $C_1 > C C_2$ with the same constant $C$.

Therefore, $T_* \leq T$. This concludes the proof of the Theorem 7.1.

Remark. The Theorem 7.1 implies that $T$ is sufficiently large:

$$T \geq \frac{1}{a}\left(\rho \log \varepsilon^{-1} - log\frac{C_1}{C_0}\right).$$

We now repeat the argument that leads to (7.12) to get the lower bound on $\|\sigma(\tau)\|_{L^Q(\mathbb{R}^2)}$.

Theorem 7.2. Given that $C_0 \in (0,1)$ is sufficiently small, the following inequality holds

$$\|\sigma(\tau) - e^{\tau L}\sigma_0\|_{L^Q(\mathbb{R}^2)} \leq CC_0C_2\varepsilon^\rho e^{a\tau}, \ \ \forall \tau \in [0, T_*].$$  (7.26)

Proof. From (7.1), using the variation of parameters

$$\sigma(\tau) = e^{\tau L}\sigma_0 + \int_0^\tau e^{(\tau-\xi)L}\left(-\big(w(\xi),\partial_y\big)\sigma(\xi)\right)d\xi \, ,$$  (7.27)

where the integral in the right side of (7.27) is understood in the sense of Bochner.

From (7.21), (7.25)

$$\|\big(w(\xi),\partial_y\big)\sigma(\xi)\|_{L^Q(\mathbb{R}^2)} \leq \||y|^{-1}w(\xi)\|_{L^\infty(\mathbb{R}^2)} \, \| \, |y|\partial_y\sigma(\xi)\|_{L^Q(\mathbb{R}^2)}$$

$$\leq \|\partial_y w(\xi)\|_{L^\infty(\mathbb{R}^2)} \, \| \, |y|\partial_y\sigma(\xi)\|_{L^Q(\mathbb{R}^2)}$$

$$\leq \|\partial_y \sigma(\xi)\|_{L^{2,1}(\mathbb{R}^2)} \, \| \, |y|\partial_y\sigma(\xi)\|_{L^Q(\mathbb{R}^2)}$$

$$\leq CC_2^2\varepsilon^{2\rho}e^{2a\xi} \, , \ \ \forall \xi \in [0, T_*].$$  (7.28)

Combining (7.27), (7.28) we get

$$\|\sigma(\tau) - e^{\tau L}\sigma_0\|_{L^Q(\mathbb{R}^2)} \leq CC_2^2\varepsilon^{2\rho}\int_0^\tau e^{a(\tau-\xi)} \, e^{2a\xi} \, d\xi$$

$$\leq CC_2^2\varepsilon^{2\rho}e^{2a\tau}$$

$$\leq CC_2^2C_1^{-1}(C_1\varepsilon^\rho e^{a\tau})\varepsilon^\rho e^{a\tau}$$

$$\leq CC_0C_2\varepsilon^\rho e^{a\tau} \, , \ \ \forall \tau \in [0, T_*].$$  (7.29)

This concludes the proof of Theorem 7.2. QED

We note in addition to (7.26) that the first term in the right side of (7.27) satisfies the following two inequalities for an appropriate $C \geq 1$ :

$$C^{-1}C_3\varepsilon^\rho e^{a\tau} \leq \|e^{\tau L}\sigma_0\|_{L^Q(\mathbb{R}^2)} \leq CC_3\varepsilon^\rho e^{a\tau}, \quad \forall \tau \in [0, T_*], \tag{7.30}$$

where

$$C_3 = \|Re\,\eta\|_{L^Q(\mathbb{R}^2)} \tag{7.31}$$

For a sufficiently small (but finite, i.e., independent of $\varepsilon$) $C_0 \in (0,1)$,

$$\|\sigma(\tau) - e^{\tau L}\sigma_0\|_{L^Q(\mathbb{R}^2)} \leq CC_0 \|e^{\tau L}\sigma_0\|_{L^Q(\mathbb{R}^2)}, \qquad \forall \tau \in [0, T_*], \tag{7.32}$$

where

$$T_* = \frac{1}{a}\Big(\rho\log\varepsilon^{-1} - log\frac{C_1}{C_0}\Big).$$

# § 8. Modified $V$

Let for a fixed $\gamma > 0$ that is small enough we choose $\rho = \gamma^{-1}a$.

$$V_1(y, \tau) = \chi(\varepsilon e^{\gamma\tau}|y|)\,V(y), \tag{8.1}$$

where $\chi$ is chosen as in (6.1).

Let $\sigma_1(y, \tau)$ be a solution to the following problem:

$$\begin{cases} \partial_\tau \sigma_1(y, \tau) = -\big(V_1(y, \tau), \partial_y\big)\sigma_1(y, \tau) - \big(w_1(y, \tau), \partial_y\big)curl_y[V_1(y, \tau)] \\ \qquad\qquad - (w_1, \partial_y)\sigma_1 + \gamma\big(\alpha\sigma_1 + (y, \partial_y)\sigma_1\big); \\ \qquad w_1(y, \tau) = (\mathcal{K}(.) * \sigma_1(.,\tau))(y), \qquad div_y w_1 = 0; \\ \qquad \sigma(y, 0) = \sigma_0(y) = Re\,(\varepsilon^{\rho + i\zeta}\,\eta(y)), \qquad y \in \mathbb{R}^2. \end{cases} \tag{8.2}$$

We proceed as in §7. We start by choosing a constant $C_2 > 0$ that depends only on $V$ so that

$$\begin{cases} \|\eta\|_{L^Q(\mathbb{R}^2)} \leq C_2; \\ \||y|\partial_y\eta\|_{L^Q(\mathbb{R}^2)} \leq C_2; \\ \|\partial_y\eta\|_{L^{2,1}(\mathbb{R}^2)} \leq C_2. \end{cases} \tag{8.3}$$

Next, we choose another constant $C_1 > CC_2$, where $C > 1$ is sufficiently large to be specified later.

Let $T \in (0, \infty]$ be *the maximal time*, so that the solution $\sigma_1(y, \tau)$ of the problem (8.2) satisfies the following three inequalities:

$$\begin{cases} \|\sigma_1(\tau)\|_{L^Q(\mathbb{R}^2)} \leq C_1 \varepsilon^\rho e^{a\tau}, & \tau \in [0, T]; \\ \||y|\partial_y \sigma_1(\tau)\|_{L^Q(\mathbb{R}^2)} \leq C_1 \varepsilon^\rho e^{a\tau}, & \tau \in [0, T]; \\ \|\partial_y \sigma_1(\tau)\|_{L^{2,1}(\mathbb{R}^2)} \leq C_1 \varepsilon^\rho e^{a\tau}, & \tau \in [0, T]. \end{cases} \tag{8.4}$$

As above in §7

$$Q \in \left(2, \frac{2}{\alpha}\right). \tag{8.5}$$

We define $T_*$ so that

$$C_1 \varepsilon^\rho e^{a T_*} = C_0 \ll 1, \tag{8.6}$$

where $C_0$ is a fixed positive constant to be chosen later.

Theorem 8.1. The following inequality holds

$$T \geq T_*. \tag{8.7}$$

Proof. We assume to the contrary that $T_* > T$.

First, we estimate $\|\sigma_1(\tau)\|_{L^Q(\mathbb{R}^2)}$ for $\tau \in [0, T]$.

From (8.2)

$$\partial_\tau \sigma_1(y, \tau) = -(V, \partial_y)\sigma_1 - (w_1, \partial_y) curl\, V$$

$$\gamma(\alpha \sigma_1 + (y, \partial_y)\sigma_1)$$

$$+((V - V_1), \partial_y)\sigma_1 + (w_1, \partial_y) curl\, (V - V_1)$$

$$-(w_1, \partial_y)\sigma_1$$

$$= L\sigma_1 + F(y, \tau). \tag{8.8}$$

Using variation of parameters formula, we get

$$\sigma_1(\tau) = e^{\tau L}\sigma_0 + \int_0^\tau e^{(\tau - \xi)L} F(\xi)\, d\xi, \quad \forall \tau \in [0, T]. \tag{8.9}$$

We need the estimate for $\|\sigma_1(\tau)\|_{L^Q(\mathbb{R}^2)}$.

From (8.4)

$$\|(w_1(\tau), \partial_y)\sigma_1(\tau)\|_{L^Q(\mathbb{R}^2)} \leq \|\frac{w_1(\tau)}{|y|}\|_{L^\infty(\mathbb{R}^2)} \quad \||y|\partial_y \sigma_1(\tau)\|_{L^Q(\mathbb{R}^2)}$$

$$\leq C\|\partial_y w_1(\tau)\|_{L^\infty(\mathbb{R}^2)}\||y|\partial_y\,\sigma_1(\tau)\|_{L^Q(\mathbb{R}^2)}$$

$$\leq C\|\partial_y\,\sigma_1(\tau)\|_{L^{2,1}(\mathbb{R}^2)}\||y|\partial_y\,\sigma_1(\tau)\|_{L^Q(\mathbb{R}^2)}$$

$$\leq CC_1^2\varepsilon^{2\rho}e^{2a\tau}\quad,\qquad \forall\tau\in[0,T]\,;\tag{8.10}$$

$$\|(V-V_1,\partial_y)\sigma_1(\tau)\|_{L^Q(\mathbb{R}^2)}\leq$$

$$\leq \||y|^{-1}(V-V_1)\|_{L^\infty(\mathbb{R}^2)}\;\;\||y|\partial_y\,\sigma_1(\tau)\|_{L^Q(\mathbb{R}^2)}$$

$$= \||y|^{-1}(1-\tilde{\chi}\;(C_1^{-1}\varepsilon e^{\gamma\tau}|y|)V(y)\|_{L^\infty(\mathbb{R}^2)}\;\||y|\partial_y\,\sigma_1(\tau)\|_{L^Q(\mathbb{R}^2)}$$

$$\leq CC_1C_4^{-\alpha}\varepsilon^\alpha e^{\alpha\gamma\tau}\varepsilon^\rho e^{a\tau},\;\forall\tau\in[0,T]\,;\tag{8.11}$$

$$\|(w_1,\partial_y)curl_y\,(V(y)-V_1(y,\tau))\|_{L^Q(\mathbb{R}^2)}$$

$$\leq C\|\frac{w_1(\tau)}{|y|}\|_{L^Q(\mathbb{R}^2)}\||y|curl_y\,(V(y)-V_1(y,\tau))\|_{L^\infty(\mathbb{R}^2)}$$

$$\leq\;C\;\|\sigma_1(\tau)\|_{L^Q(\mathbb{R}^2)}\||y|\partial_y curl_y\,[(1-\tilde{\chi}\;(C_4^{-1}\varepsilon e^{\gamma\tau}|y|)V(y)]\|_{L^\infty(\mathbb{R}^2)}$$

$$\leq CC_1C_4^{-\alpha}\varepsilon^\alpha e^{\alpha\gamma\tau}\varepsilon^\rho e^{a\tau},\;\forall\tau\in[0,T].\tag{8.12}$$

From (8.9) - (8.12)

$$\|\sigma_1(\tau)\|_{L^Q(\mathbb{R}^2)}\leq C_2\varepsilon^\rho e^{a\tau}+CC_1C_4^{-\alpha}\int_0^\tau e^{a(\tau-\xi)}\varepsilon^{\alpha+\rho}\,e^{(\alpha\gamma+a)\xi}\,d\xi$$

$$+CC_1^2\;\int_0^\tau e^{a(\tau-\xi)}\varepsilon^{2\rho}\,e^{2a\xi}\,d\xi$$

$$\leq C_2\varepsilon^\rho e^{a\tau}+CC_1C_4^{-\alpha}\varepsilon^{\alpha+\rho}e^{(\alpha\gamma+a)\tau}+CC_1^2e^{2\rho}e^{2a\tau}$$

$$\leq C_2\varepsilon^\rho e^{a\tau}+C_1^{1-\frac{\alpha}{\rho}}C_4^{-\alpha}C_0^{\frac{\alpha}{\rho}}\varepsilon^\rho e^{a\tau}+CC_0C_1\varepsilon^\rho e^{a\tau}$$

$$\leq 2C_2\;\varepsilon^\rho e^{a\tau},\qquad \forall\tau\in[0,T],\tag{8.13}$$

for $C_0\in(0,1)$ sufficiently small.

We now differentiate (8.2) to get

$$\partial_\tau \partial_k \sigma_1(y,\tau) = -\big(V_1(y,\tau),\partial_y\big)\partial_k\sigma_1(y,\tau) - \big(\partial_k V_1(y,\tau),\partial_y\big)\sigma_1$$

$$+\gamma\big((\alpha+1)\partial_k\sigma_1 + (y,\partial_y)\partial_k\sigma_1\big)$$

$$-\big(\partial_k w_1(y,\tau),\partial_y\big)curl_y[V_1(y,\tau)] - \big(w_1(y,\tau),\partial_y\big)\partial_k curl_y[V_1(y,\tau)]$$

$$-\big(w_1,\partial_y\big)\partial_k\sigma_1 - \big(\partial_k w_1,\partial_y\big)\sigma_1;$$

$$\partial_k\sigma_1(y,0) = \partial_k\sigma_0(y), \ \ k=1,2\,;$$

$$w_1(y,\tau) = (\mathcal{K}(.)*\sigma_1(.\,,\tau)(y)\,, \quad div_y w_1 = 0\,. \tag{8.14}$$

**Lemma 8.1.** Let

$$\begin{cases} \dfrac{dx}{d\tau} = V_1(x,\tau) + w_1(x,\tau) - \gamma x; \\ \qquad x(0) = y \in \mathbb{R}^2. \end{cases}$$

Then, for any $\tau \in [0,T]$ and given that $C_0$ is sufficiently small, the following inequality holds:

$$\int_0^\tau |\frac{1}{2}\Big(\big(\partial_y V_1(x(\xi),\xi)\big) + \big(\partial_y V_1(x(\xi),\xi)\big)^*\Big)|d\xi \le C < \infty\,,$$

where the constant $C > 0$ in (7.13) is uniform with respect to $y \in \mathbb{R}^2$. Here $(\partial_y V_1)^*$ stands for the transposition of the matrix $\partial_y V_1$.

Proof. We have

$$\frac{d}{d\tau}|x(\tau)| = \Big(\frac{w_1(x(\tau),\tau)}{|x(\tau)|},x(\tau)\Big) - \gamma|x(\tau)|; \tag{8.15}$$

$$\|\frac{w_1(.,\tau)}{|.|}\|_{L^\infty(\mathbb{R}^2)} \le C\|\partial_y w_1(\tau)\|_{L^\infty(\mathbb{R}^2)}$$

$$\le C\|\partial_y\sigma_1(\tau)\|_{L^{1,2}(\mathbb{R}^2)} \le CC_1\varepsilon^\rho e^{\alpha\tau}, \ \ \forall \tau \in [0,T]. \tag{8.16}$$

From (7.14), (7.15) we derive the following two inequalities

$$|y|\exp\{-\gamma\tau - CC_1\varepsilon^\rho e^{\alpha\tau}\} \le |x(\tau)| \le |y|\exp\{-\gamma\tau + CC_1\varepsilon^\rho e^{\alpha\tau}\}, \quad \forall \tau \in [0,T] \tag{8.17}$$

Also, $V_1(y) = \Omega_1(|y|, \tau)y^\perp$, where $\Omega_1(|y|, \tau) = \chi(\varepsilon e^{\gamma\tau}|y|)\Omega(|y|)$ so

$$\frac{1}{2}\left(\left(\partial_y V_1(y, \tau)\right) + \left(\partial_y V_1(y, \tau)\right)^*\right) = \frac{1}{|y|}\,\partial_{|y|}\Omega_1(|y|, \tau)\begin{pmatrix} -y_1 y_2 & \frac{1}{2}\left(y_1^2 - y_2^2\right) \\ \frac{1}{2}\left(y_1^2 - y_2^2\right) & y_1 y_2 \end{pmatrix},$$

therefore,

$$\left|\frac{1}{2}\left(\left(\partial_y V_1(y, \tau)\right) + \left(\partial_y V_1(y, \tau)\right)^*\right)\right| \le C|y|^2(1 + |y|)^{-2-\alpha}. \tag{8.18}$$

Here we assume $C_0 \in (0,1)$ is sufficiently small   and $\tau \in [0, T]$.

If $C_0$ is chosen to be small enough, so that $CC_1 \varepsilon^\rho e^{\alpha\tau} \le CC_0 \le 1$ in (8.17), then the integral in the statement of the lemma is

$$\le C \int_0^\tau \left|y e^{-\gamma\xi}\right|^2 \left(1 + e^{-\gamma\xi}|y|\right)^{-2-\alpha} d\xi$$

$$= \frac{1}{\gamma} C \int_{e^{-\gamma\tau}|y|}^{|y|} r \left(1 + r\right)^{-2-\alpha} dr$$

$$\le \frac{1}{\gamma} C \int_0^\infty r \left(1 + r\right)^{-2-\alpha} dr.$$

Since the last integral converges, this completes the proof of Lemma 8.1. QED.

We now estimate $\|\partial_y \sigma(\tau)\|_{L^{2,1}(\mathbb{R}^2)}$ using the system (8.14). We apply Lemma 7.3 with

$$W_k(y, \tau) = \partial_k \sigma_1(y, \tau);\ k = 1,2\ ;\ n = 2, \qquad N = 2;$$

$$Z(y, \tau) = V_1(y, \tau) + w_1(y, \tau) - \gamma y;$$

$$\mu = -\gamma(\alpha + 1), \nu = -2\gamma;$$

$$A_{kl}(y, \tau) = \partial_k V_{1l}(y, \tau) + \partial_k w_{1l}(y, \tau);\ k, l = 1,2;$$

$$F_k(y, \tau) = -\left(\partial_k w_1, \partial_y\right) curl\ V_1(y, \tau) - \left(w_1, \partial_y\right)\partial_k curl\ V_1(y, \tau);\ k = 1,2;$$

$$g(\tau) = 0\ , \qquad \forall \tau \in [0, T].$$

From Lemma 8.1

$$\int_0^\tau \frac{1}{2}|(A+A^*)(x(\xi),\xi)|d\xi \quad \leq \quad \int_0^\tau \left|\frac{1}{2}\Big(\big(\partial_y V_1\big)+\big(\partial_y V_1\big)^*\Big)\right|(x(\xi),\xi)\,d\xi$$

$$+\int_0^\tau \left|\frac{1}{2}\Big(\big(\partial_y w_1\big)+\big(\partial_y w_1\big)^*\Big)\right|(x(\xi),\xi)\,d\xi$$

$$\leq C+\int_0^\tau \|\partial_y w_1(\xi)\|_{L^\infty(\mathbb{R}^2)}\,d\xi$$

$$\leq C+\int_0^\tau \|\partial_y \sigma_1(\xi)\|_{L^{1,2}(\mathbb{R}^2)}\,d\xi$$

$$\leq C+CC_1\int_0^\tau \varepsilon^\rho\, e^{a\xi}\,d\xi$$

$$\leq C+CC_1\varepsilon^\rho e^{a\tau} \leq C+CC_0 \quad,\quad \forall \tau \in [0,T]\,. \tag{8.19}$$

To estimate $\|F(\tau)\|_{L^{2,1}(\mathbb{R}^2)}$ we proceed as in Section 7:

$$\|\partial_k w_1\|_{L^Q(\mathbb{R}^2)}\leq C\|\sigma_1\|_{L^Q(\mathbb{R}^2)}\,,\quad k=1,2\,;$$

$$\|\frac{1}{|y|}w_1\|_{L^Q(\mathbb{R}^2)}\leq C\|\sigma_1\|_{L^Q(\mathbb{R}^2)},\quad Q\in\left(2,\frac{2}{\alpha}\right)\quad;$$

$$\begin{cases} \|\big(\partial_k w_1,\partial_y\big)curl\,V_1(y)\|_{L^{2,1}(\mathbb{R}^2)} \;\leq\; C\|\partial_k w_1\|_{L^Q(\mathbb{R}^2)}\|\partial_y curl\,V_1\|_{L^{\frac{2Q}{Q-2},\frac{Q}{Q-1}}(\mathbb{R}^2)} \leq C\|\sigma_1\|_{L^Q(\mathbb{R}^2)} \\[2mm] \|(w_1,\partial_y)\partial_k curl\,V_1\|_{L^{2,1}(\mathbb{R}^2)}\leq C\|\frac{1}{|y|}w_1\|_{L^Q(\mathbb{R}^2)}\||y|\partial_y\partial_k curl\,V_1\|_{L^{\frac{2Q}{Q-2},\frac{Q}{Q-1}}(\mathbb{R}^2)} \leq C\|\sigma_1\|_{L^Q(\mathbb{R}^2)} \end{cases} \tag{8.20}$$

From (8.20), (8.13)

$$\|F(\tau)\|_{L^{2,1}(\mathbb{R}^2)}\leq CC_2\varepsilon^\rho e^{a\tau},\;\;\forall\tau\in[0,T]. \tag{8.21}$$

From (8.14), Lemma 7.3 and (8.19) --(8.21) we derive:

$$\|\partial_y\sigma_1(\tau)\|_{L^{2,1}(\mathbb{R}^2)}\leq CC_2\varepsilon^\rho e^{\alpha\gamma\tau}+CC_2\int_0^\tau e^{\alpha\gamma(\tau-\xi)}\,\varepsilon^\rho e^{a\xi}\,d\xi\leq CC_2\varepsilon^\rho e^{a\tau},\quad \forall\tau\in[0,T], \tag{8.22}$$

provided $\gamma<\frac{1}{\alpha}a.$

We turn now to the estimate of $\||y|\partial_y\sigma_1(\tau)\|_{L^Q(\mathbb{R}^2)},\;\tau\in[0,T].$

From (8.14) we get the system of equations (summation with respect to the indexes repeated twice is assumed):

$$\partial_\tau(|y|\partial_k\sigma_1(y,\tau)) = -(V_1,\partial_y)(|y|\partial_k\sigma_1) - \partial_k V_{1j}|y|\partial_j\sigma_1$$

$$+\gamma(\alpha|y|\partial_k\sigma_1 + (y,\partial_y)(|y|\partial_k\sigma_1))$$

$$-\partial_k w_{1j}|y|\partial_j\,curl\,V_1 - w_{1j}|y|\partial_j\partial_k\,curl\,V_1$$

$$-(w_1,\partial_y)(|y|\partial_k\sigma_1) + |y|^{-2}(w_1,y)(|y|\partial_k\sigma_1)$$

$$-\partial_k w_{1j}(|y|\partial_j\sigma_1). \tag{8.23}$$

We use Lemma 7.2 with $q = Q$,

$$W_k(y,\tau) = |y|\,\partial_k\sigma_1(y,\tau); \quad k = 1,2 \,; N = 2, n = 2;$$

$$Z(y,\tau) = V_1(y,\tau) + w_1(y,\tau) - \gamma y;$$

$$\mu = -\alpha\gamma, \nu = -2\gamma;$$

$$A_{kj}(y,\tau) = \partial_k V_{1j}(y,\tau) + \partial_k w_{1j}(y,\tau) - |y|^{-2}(w_1,y)\delta_{kj}; \quad k,j = 1,2;$$

$$g(\tau) = 0, \qquad \forall\tau \in [0,T];$$

$$F_k(y,\tau) = -\partial_k w_{1j}|y|\partial_j\,curl_y\,V_1(y,\tau) - w_{1j}|y|\partial_j\partial_k\,curl_y\,V_1(y,\tau), k = 1,2.$$

Here $\delta_{kj}$ is Kronecker delta.

We first estimate $\|F(\tau)\|_{L^Q(\mathbb{R}^2)}$. We have

$$\|F(\tau)\|_{L^Q(\mathbb{R}^2)} \le \|\partial_y w_1(\tau)\|_{L^Q(\mathbb{R}^2)} \; \||y|\partial_y\,curl_y\,V_1(\tau)\|_{L^\infty(\mathbb{R}^2)}$$

$$+\||y|^{-1}w_1(\tau)\|_{L^Q(\mathbb{R}^2)}\||y|^2\partial_y^2\,curl_y\,V_1(\tau)\|_{L^\infty(\mathbb{R}^2)}$$

$$\le C\|\sigma_1(\tau)\|_{L^Q(\mathbb{R}^2)}, \tag{8.24}$$

since $Q \in \left(2, \frac{2}{\alpha}\right) \subset (2,\infty)$.

The integral

$$\int_0^\tau \frac{1}{2}|(A + A^*)(x(\xi),\xi)|d\xi$$

contains one additional term as compared to the previous inequality (8.19). This term is bounded above by

$$\int_0^\tau \||y|^{-1}w_1(\xi)\|_{L^\infty(\mathbb{R}^2)}\,d\xi$$

$$\leq C\int_0^\tau \|\partial_y w_1(\xi)\|_{L^\infty(\mathbb{R}^2)}\,d\xi$$

$$\leq C\int_0^\tau \|\partial_y \sigma_1(\xi)\|_{L^{2,1}(\mathbb{R}^2)}\,d\xi$$

$$\leq C\int_0^\tau C_2\varepsilon^\rho e^{a\xi}\,d\xi$$

$$\leq CC_2\varepsilon^\rho e^{a\tau} \leq CC_0C_2C_1^{-1} \leq CC_0, \ \ \forall \tau \in [0,T].$$

Therefore,

$$\int_0^\tau \frac{1}{2}|(A+A^*)(x(\xi),\xi)|d\xi \ \ \leq C < \infty, \qquad \forall \tau \in [0,T]. \tag{8.25}$$

Using (8.23) and (8.24) -(8.25), (8.13) we arrive at the inequalities

$$\||y|\partial_y\sigma_1(\tau)\|_{L^Q(\mathbb{R}^2)} \leq$$

$$\leq CC_2\varepsilon^\rho e^{\gamma\left(\alpha-\frac{2}{Q}\right)\tau} + C\int_0^\tau e^{\gamma\left(\alpha-\frac{2}{Q}\right)(\tau-\xi)}\|\sigma_1(\xi)\|_{L^Q(\mathbb{R}^2)}d\xi$$

$$\leq CC_2\varepsilon^\rho e^{\gamma\left(\alpha-\frac{2}{Q}\right)\tau} + CC_2\varepsilon^\rho e^{a\tau}$$

$$\leq CC_2\varepsilon^\rho e^{a\tau}, \quad \forall \tau \in [0,T]. \tag{8.26}$$

Here we used the inequalities $\ \gamma > 0, 2 < Q < \frac{2}{\alpha}$.

From (8.4) and the definition of $T$, at least one of the three inequalities in (8.4) should become an equality for $\tau = T$. Because of the inequalities (8.13), (8.22), (8.26),

$$C_1 \leq CC_2.$$

But this is impossible, provided that $C_1$ had been chosen to be large enough ( $C_1 > CC_2$ with a sufficiently large constant $C > 1$ ) to start with.

This concludes the proof of the Theorem 8.1. QED.

Theorem 8.2. Let as above,

$$L\sigma \equiv -(V(y), \partial_y)\sigma - (w, \partial_y)\, curl\, V(y) + \gamma(\alpha\sigma + (y, \partial_y)\sigma);$$

$$w = \mathcal{K} *_y \sigma,$$

Where $\gamma > 0$ is fixed and small enough. Given that $C_0 > 0$ is sufficiently small, the following inequality holds:

$$\|\sigma_1(\tau) - e^{\tau L}\sigma_0\|_{L^Q(\mathbb{R}^2)} \le C C_0^{\frac{\alpha}{\rho}} C_1^{1-\frac{\alpha}{\rho}} \varepsilon^\rho e^{a\tau}, \ \ \forall \tau \in [0, T_*].$$

Proof. From (8.8)

$$\partial_\tau \sigma_1 = -(V, \partial_y)\sigma_1 - (w_1, \partial_y)curl\, V$$

$$+\gamma(\alpha\sigma_1 + (y, \partial_y)\sigma_1)$$

$$+((V - V_1), \partial_y)\sigma_1 + (w_1, \partial_y)curl\, (V - V_1)$$

$$-(w_1, \partial_y)\sigma_1$$

$$= L\sigma_1 + F(y, \tau). \tag{8.27}$$

We first estimate $\|F(\tau)\|_{L^Q(\mathbb{R}^2)}, \ \forall \tau \in [0, T_*]$.

From (8.4) for $\tau \in [0, T_*]$ using the Theorem 8.1,

$$\|(w_1(\tau), \partial_y)\sigma_1(\tau)\|_{L^Q(\mathbb{R}^2)} \le C C_1^2 \varepsilon^{2\rho} e^{2a\tau} \quad , \tag{8.28}$$

as in (8.10);

$$\|(V - V_1, \partial_y)\sigma_1(\tau)\|_{L^Q(\mathbb{R}^2)} \le C C_1 C_4^{-\alpha} \varepsilon^\alpha e^{\alpha\gamma\tau} \varepsilon^\rho e^{a\tau}, \tag{8.29}$$

as in (8.11). Also,

$$\|(w_1(\tau), \partial_y)curl_y\, (V(y) - V_1(y, \tau))\|_{L^Q(\mathbb{R}^2)}$$

$$\le C C_1 C_4^{-\alpha} \varepsilon^\alpha e^{\alpha\gamma\tau} \varepsilon^\rho e^{a\tau}, \tag{8.30}$$

as in (8.12).

From (8.27), (8.28) --(8.30)

$$\|\sigma_1(\tau) - e^{\tau L}\sigma_0\|_{L^Q(\mathbb{R}^2)}$$

$$\leq \| \int_0^\tau e^{(\tau-\xi)L} \, F(\xi) \, d\xi \|_{L^Q(\mathbb{R}^2)}$$

$$\leq C \int_0^\tau e^{a(\tau-\xi)} (C_1^2 \varepsilon^{2\rho} e^{2a\xi} + C_1 C_4^{-\alpha} \varepsilon^\alpha e^{\alpha\gamma\xi} \varepsilon^\rho e^{a\xi}) d\xi$$

$$\leq CC_1^2 \varepsilon^{2\rho} e^{2a\tau} + C_1 C_4^{-\alpha} \varepsilon^{\alpha+\rho} e^{(\gamma\alpha+a)\tau}, \quad \forall \tau \in [0, T_*]. \tag{8.31}$$

The inequality (8.31) yields

$$\|\sigma_1(\tau) - e^{\tau L}\sigma_0\|_{L^Q(\mathbb{R}^2)}$$

$$\leq CC_0 C_1 \varepsilon^\rho e^{a\tau} + CC_0^{\frac{\alpha}{\rho}} C_1^{1-\frac{\alpha}{\rho}} \varepsilon^\rho e^{a\tau}$$

$$\leq CC_0^{\frac{\alpha}{\rho}} C_1^{1-\frac{\alpha}{\rho}} \varepsilon^\rho e^{a\tau}, \quad \forall \tau \in [0, T_*]. \tag{8.32}$$

This concludes the prof of Theorem 8.2. QED.

Remark. As above,

$$C^{-1} C_3 \varepsilon^\rho e^{a\tau} \leq \|e^{\tau L}\sigma_0\|_{L^Q(\mathbb{R}^2)} \leq CC_3 \varepsilon^\rho e^{a\tau}, \qquad \forall \tau \in [0, T_*], \tag{8.33}$$

where

$$C_3 = \|Re \, \eta\|_{L^Q(\mathbb{R}^2)}$$

Therefore, for a sufficiently small (but finite, i.e., independent of $\varepsilon$) $C_0 \in (0,1)$,

$$\|\sigma_1(\tau) - e^{\tau L}\sigma_0\|_{L^Q(\mathbb{R}^2)} \leq CC_0^{\frac{\alpha}{\rho}} C_1^{1-\frac{\alpha}{\rho}} \|e^{\tau L}\sigma_0\|_{L^Q(\mathbb{R}^2)}, \qquad \forall \tau \in [0, T_*], \tag{8.34}$$

where

$$T_* = \frac{1}{a}\left(\rho \log \varepsilon^{-1} - log \frac{C_1}{C_0}\right).$$

# § 9. The limit $\varepsilon \to 0$

We now employ the fundamental scaling properties of the Euler equations of an ideal incompressible fluid.

Let $\omega(x,t), v(x,t) = (\mathcal{K}(.) * \omega(.,t))(x)$ be a (strong) solution of the Euler equations (4.1).

For any $\varepsilon > 0$ let

$$\begin{cases} \omega_\varepsilon(x,t) = \varepsilon^{-\alpha}\omega(\varepsilon^{-1}x, \varepsilon^{-\alpha}t); \\ v_\varepsilon(x,t) = \varepsilon^{1-\alpha}v(\varepsilon^{-1}x, \varepsilon^{-\alpha}t); \\ Z_\varepsilon(x,t) = \varepsilon^{-2\alpha}Z(\varepsilon^{-1}x, \varepsilon^{-\alpha}t); \end{cases} \tag{9.1}$$

$$\omega_{0\varepsilon}(x) = \varepsilon^{-\alpha}\omega_0(\varepsilon^{-1}x). \tag{9.2}$$

With this rescaling (9.1) satisfies the Euler equations

$$\begin{cases} \partial_t\omega_\varepsilon(x,t) = -(v_\varepsilon(x,t),\partial_x)\omega_\varepsilon(x,t) + Z_\varepsilon(x,t); \ x \in \mathbb{R}^2, \ t \in [0,\infty), \\ \omega_\varepsilon(x,t) = curl_x \ v_\varepsilon(x,t); \ div_x v_\varepsilon(x,t) = 0, \\ \omega_\varepsilon(x,0) = \omega_{0\varepsilon}(x). \end{cases} \tag{9.3}$$

We now choose $\omega(x,t)$ as in (6.3).

The external force $Z(x,t)$ is given by (6.8) under the condition (6.7). From (6.3), (6.9), (6.12), (6.19) we have

$$\omega_0(x) = Re\left(\varepsilon^{\rho+i\zeta}\eta(x)\right) + curl_x\left[\chi(\varepsilon|x|)V(x)\right]. \tag{9.4}$$

Therefore,

$$\omega_{0\varepsilon}(x) = Re(\varepsilon^{\rho+i\zeta-\alpha}\eta(\varepsilon^{-1}x))$$

$$+\chi(|x|)\,\varepsilon^{-\alpha}(curl\,V)(\varepsilon^{-1}x) + |x|\,\chi'(|x|)\varepsilon^{-\alpha}\Omega(\varepsilon^{-1}x)$$

$$= Re(\varepsilon^{\rho+i\zeta-\alpha}\eta(\varepsilon^{-1}x)) + curl_x\left[\chi(|x|)\,\varepsilon^{1-\alpha}\,V(\varepsilon^{-1}x)\right]. \tag{9.5}$$

Since $(\omega_\varepsilon(x,t), v_\varepsilon(x,t))$ is a strong solution as in (4.1), the following integral identity holds for every

$$\psi(x,t) \in C_0^\infty\left([0,\infty); S(\mathbb{R}^2)\right):$$

$$\int \omega_{0\varepsilon}(x)\,\psi(x,0)\,dx + \int\limits_0^\infty \int \omega_\varepsilon(x,t)\,\partial_t\,\psi(x,t)\,dx\,dt$$

$$\int_0^\infty \int \omega_\varepsilon(x,t)\ (v_\varepsilon(x,t),\partial_x)\psi(x,t)\,dx\,dt + \int_0^\infty \int Z_\varepsilon(x,t)\,\psi(x,t)\,dx\,dt = 0. \tag{9.6}$$

Remark. $\varepsilon = 1$ in (9.6) corresponds to the weak formulation for the solution $(\omega(x,t),v(x,t))$. We study (9.3). The basic observation is the following.

Theorem 9.1. The following estimate holds for any $T > 0,\ \ Q \in [1,\frac{2}{\alpha})$

$$\int_0^T \| Z_\varepsilon(t)\|_{L^Q(\mathbb{R}^2)}\ dt \le C < \infty\,,$$

where the constant $C$ is uniform with respect to $\varepsilon \in (0,1]$.

Proof. Using (6.5),

$$Z(x,t) = \gamma\,\mathcal{R}(t)^{-2\alpha}(\chi(\varepsilon|x|)W(\mathcal{R}(t)^{-1}x) + \varepsilon|x|\chi'(\varepsilon|x|)B(\mathcal{R}(t)^{-1}x))$$

$$= \gamma(1+\alpha\gamma t)^{-2}(\chi(\varepsilon|x|)W((1+\alpha\gamma t)^{-\frac{1}{\alpha}}\,x) + \ \varepsilon|x|\chi'(\varepsilon|x|)B((1+\alpha\gamma t)^{-\frac{1}{\alpha}}\,x)).$$

Therefore,

$$Z_\varepsilon(x,t) = \gamma\varepsilon^{-2\alpha}(1+\alpha\gamma\varepsilon^{-\alpha}t)^{-2}(\chi(|x|)W((1+\alpha\gamma\varepsilon^{-\alpha}t)^{-\frac{1}{\alpha}}\varepsilon^{-1}\,x)$$

$$+ |x|\chi'(|x|)B((1+\alpha\gamma\varepsilon^{-\alpha}t)^{-\frac{1}{\alpha}}\varepsilon^{-1}\,x))$$

$$= \gamma(\varepsilon^\alpha+\alpha\gamma t)^{-2}(\chi(|x|)W((\varepsilon^\alpha+\alpha\gamma t)^{-\frac{1}{\alpha}}\,x)$$

$$+|x|\chi'(|x|)B((\varepsilon^\alpha+\alpha\gamma t)^{-\frac{1}{\alpha}}\,x))\,. \tag{9.8}$$

The function $x\ \rightarrow\ |x|\,\chi'(|x|)B((\varepsilon^\alpha+\alpha\gamma t)^{-\frac{1}{\alpha}}\,x))$ is supported on $\{\,C_4 \le |x| \le 2C_4\}$ and

$$\||x|\,\chi'(|x|)B((\varepsilon^\alpha+\alpha\gamma t)^{-\frac{1}{\alpha}}\,x))\|_{L^\infty(\mathbb{R}^2)} \le C(\varepsilon^\alpha+\alpha\gamma t)^{\frac{2}{\alpha}}.$$

Therefore, for any $Q \in [1,\infty]$,

$$\| |x| \, \chi'(|x|) B((\varepsilon^{\alpha} + \alpha \gamma t)^{-\frac{1}{\alpha}} \, x)) \|_{L^{Q}(\mathbb{R}^2)} \leq C(\varepsilon^{\alpha} + \alpha \gamma t)^{\frac{2}{\alpha}} . \tag{9.9}$$

Here the constant $C$ is uniform with respect to $Q \in [1, \infty]$.

Also, $W(y) \in C_0^{\infty}(\mathbb{R}^2)$, thus

$$\| (\varepsilon^{\alpha} + \alpha \gamma t)^{-2} (\chi(|x|) W((\varepsilon^{\alpha} + \alpha \gamma t)^{-\frac{1}{\alpha}} \, x) \|_{L^{Q}(\mathbb{R}^2)} \leq C(\varepsilon^{\alpha} + \alpha \gamma t)^{\frac{2}{\alpha Q} - 2}, \tag{9.10}$$

with the constant $C$ uniform with respect to $Q \in [1, \infty]$.

Therefore,

$$\int_0^T \| Z_{\varepsilon}(t) \|_{L^{Q}(\mathbb{R}^2)} \, dt$$

$$\leq C \int_0^T (\varepsilon^{\alpha} + \alpha \gamma t)^{\frac{2}{\alpha} - 2} \, dt \; + \; C \int_0^T (\varepsilon^{\alpha} + \alpha \gamma t)^{\frac{2}{\alpha Q} - 2} \, dt$$

$$\leq C \int_0^T (1 + \alpha \gamma t)^{\frac{2}{\alpha} - 2} \, dt + \; C \int_0^T (\alpha \gamma t)^{\frac{2}{\alpha Q} - 2} \, dt$$

$$\leq C < \infty,$$

Since $\alpha Q < 2$.

This concludes the proof of the Theorem 9.1. QED.

Proposition 9.1. As $\varepsilon \to 0$,

$$Z_{\varepsilon}(t) \to \bar{Z}(t)$$

strongly in $L_{loc}^1([0, \infty); L^{Q}(\mathbb{R}^2))$, $\forall Q \in \left[1, \frac{2}{\alpha}\right)$ ,

where

$$\bar{Z}(x, t) = \gamma (\alpha \gamma t)^{-2} (\chi(|x|) W((\alpha \gamma t)^{-\frac{1}{\alpha}} \, x) + |x| \chi'(|x|) B((\alpha \gamma t)^{-\frac{1}{\alpha}} \, x))$$

$$= \frac{\partial}{\partial t} curl_x \left[ \chi(|x|) (\alpha \gamma t)^{-1 + \frac{1}{\alpha}} \, V \left( (\alpha \gamma t)^{-\frac{1}{\alpha}} x \right) \right].$$

Proof. The Theorem follows immediately from (9.8) --(9.10) and the Lebesgue dominated convergence theorem. QED.

Remark. The same argument yields strong convergence $Z_\varepsilon(t) \to Z_0(t)$ in

$L^r_{loc}\left([0,\infty);\ L^Q(\mathbb{R}^2)\right),\quad \forall\, Q \in \left[1,\frac{2}{\alpha}\right)\ ;\ \forall\, r \in \left[1,\frac{\alpha Q}{2(1-\alpha Q)}\right) for\ \alpha Q \in [1,2), \forall\, r \in [1,\infty) for\ \alpha Q \in (0,1).$

For the family of strong (in fact, $C^\infty$ ) solutions of the problem (4.1) the following apriori estimate holds.

Theorem 9.2. For any $Q \in [1,\infty],\quad t \geq 0$

$$\|\omega(t)\|_{L^Q(\mathbb{R}^2)} \leq \|\omega_0\|_{L^Q(\mathbb{R}^2)} + \int_0^t \|Z(s)\|_{L^Q(\mathbb{R}^2)} ds \quad . \tag{9.13}$$

Proof. Let for $s \in [0,t]\quad g(t,s)y,\ y \in \mathbb{R}^2$  be the solution to the Cauchy problem

$$\begin{cases} \dot{x} = v(x,t); \\ x(s) = y; \\ g(t,s)y = x(t). \end{cases}$$

For any $t \geq 0$ , $s \in [0,t]$ the mapping $y \to g(t,s)y$ is a volume preserving diffeomorphism $\mathbb{R}^2 \to \mathbb{R}^2$.

Solving (4.1) by the method of characteristics,

$$\omega(x,t) = \omega_0(g(t,0)^{-1}x) + \int_0^t Z(g(t,s)^{-1}x,s)\ ds\ . \tag{9.14}$$

From (9.14) using Minkowski integral inequality,

$$\|\omega(t)\|_{L^Q(\mathbb{R}^2)} \leq \|\omega_0 \circ g(t,0)^{-1}\|_{L^Q(\mathbb{R}^2)} + \int_0^t \|Z(.,s)\circ g(t,s)^{-1}\|_{L^Q(\mathbb{R}^2)}\ ds\ .$$

But any volume preserving diffeomorphism $\mathbb{R}^2 \to \mathbb{R}^2$ conserves the $L^Q$ norm of a function. Thus, the inequality (9.13) follows. QED.

Proposition 9.2. For any  $T > 0$ and for any $Q \in [1,\frac{2}{\alpha})$

$$\|\omega_\varepsilon(t)\|_{L^Q(\mathbb{R}^2)} \leq C < \infty,\qquad \forall t \in [0,T], \tag{9.15}$$

where the constant $C$ may depend on $T$  but is uniform with respect to   $\varepsilon \in (0,1]$.

Proof. From (9.3) and Theorem 9.2

$$\|\omega_\varepsilon(t)\|_{L^Q(\mathbb{R}^2)} \leq \|\omega_{0\varepsilon}\|_{L^Q(\mathbb{R}^2)} + \int_0^t \|Z_\varepsilon(s)\|_{L^Q(\mathbb{R}^2)} ds$$

$$\leq C\varepsilon^{\rho-\alpha+\frac{2}{Q}} + \int_0^t \|Z_\varepsilon(s)\|_{L^Q(\mathbb{R}^2)} ds \quad \leq C < \infty, \ \forall \ t \in [0,T].$$

This concludes the proof of the Proposition. QED.

To study passing to the limit as $\varepsilon \to 0$ in the nonlinear term in (9.6), we need the following proposition.

Proposition 9.3. Let for $x \in \mathbb{R}^n$, $t \in [0,\infty)$

$$u(x,t) \in C^\infty\big([0,\infty); S'(\mathbb{R}^n)\big) \cap L_{loc}^\infty\big([0,\infty); L^Q(\mathbb{R}^n)\big);$$

$$\partial_x u(x,t) \in L_{loc}^\infty\big([0,\infty); L^Q(\mathbb{R}^n)\big);$$

$$\partial_t u(x,t) \in L_{loc}^1\big([0,\infty); L^Q(\mathbb{R}^n)\big),$$

where $Q \in (n,\infty)$. Assume

$$\|\partial_t u(x,t)\|_{L^Q(\mathbb{R}^2)} \leq f(t),$$

where

$$f \in L_{loc}^1\big([0,\infty)\big).$$

Then, for any, $T > 0$ the function $u(x,t)$ is continuous on the set $\mathbb{R}^n \times [0,T]$ with a qualified modulus of continuity, that depends only on $f$. The following inequality holds:

$$|u(x_2,t_2) - u(x_1,t_1)| \leq C\left(|x_2 - x_1|^{1-\frac{n}{Q}} + \left(\int_{t_1}^{t_2} f(t)dt\right)^{1-\frac{n}{Q}}\right),$$

where $x_1, x_2 \in \mathbb{R}^n, |x_1 - x_2| \leq 1, \quad 0 \leq t_1 \leq t_2 \leq T, \quad \int_{t_1}^{t_2} f(t)dt \leq 1.$

Proof. Let

$$u(x,t) = \sum_{j=-1}^\infty \Delta_j u(x,t)$$

be the Littlewood-Paley decomposition of $u(x,t)$ in $x$ −variable.

The Theorem follows from two inequalities (the first one is very well known):

$$|u(x_2, t) - u(x_1, t)| \leq C|x_2 - x_1|^{1-\frac{n}{Q}} \ , \forall t \in [0, T]; \tag{9.16}$$

$$|u(x, t_2) - u(x, t_1)| \leq C \left( \int_{t_1}^{t_2} f(t) dt \right)^{1-\frac{n}{Q}} \ , \quad \forall \ x \in \mathbb{R}^n, \tag{9.17}$$

Where $x_1, x_2 \in \mathbb{R}^n, |x_1 - x_2| \leq 1, \quad 0 \leq t_1 \leq t_2 \leq T, \quad \int_{t_1}^{t_2} f(t) dt \leq 1.$

We have

$$u(x_2, t) - u(x_1, t) = \sum_{j=-1}^{N} \left( \Delta_j u(x_2, t) - \Delta_j u(x_1, t) \right) + \sum_{j=N+1}^{\infty} \left( \Delta_j u(x_2, t) - \Delta_j u(x_1, t) \right),$$

where $N$ will be chosen later. Therefore, from Bernstein's inequality,

$$|u(x_2, t) - u(x_1, t)| \leq |x_1 - x_2| \ \| \sum_{j=-1}^{N} \Delta_j \partial_x u(. \ , t) \ \|_{L^\infty(\mathbb{R}^n)}$$

$$+ 2 \sum_{j=N+1}^{\infty} \| \Delta_j u(. \ , t) \|_{L^\infty(\mathbb{R}^n)}$$

$$\leq C|x_1 - x_2| \ \| \sum_{j=-1}^{N} \Delta_j \partial_x u(. \ , t) \ \|_{L^Q(\mathbb{R}^n)} 2^{\frac{Nn}{Q}}$$

$$+ C \sum_{j=N+1}^{\infty} \| \Delta_j \partial_x u(. \ , t) \|_{L^Q(\mathbb{R}^n)} \ 2^{\frac{jn}{Q}-j}$$

$$\leq C|x_1 - x_2| \| \partial_x u \|_{L^\infty([0,T]; \ L^Q(\mathbb{R}^n))} 2^{\frac{Nn}{Q}}$$

$$+ C \| \partial_x u \|_{L^\infty([0,T]; \ L^Q(\mathbb{R}^n))} 2^{\frac{Nn}{Q}-N}.$$

Choosing $N = [-log_2|x_1 - x_2|]$ we arrive at (9.16).

To prove (9.17) we use the Littlewood-Paley decomposition of $u(x, t)$:

$$u(x, t_2) - u(x, t_1) = \sum_{j=-1}^{M} \left( \Delta_j u(x, t_2) - \Delta_j u(x, t_1) \right) + \sum_{j=M+1}^{\infty} \left( \Delta_j u(x, t_2) - \Delta_j u(x, t_1) \right),$$

where $M$ will be chosen later. Therefore,

$$|u(x, t_2) - u(x, t_1)| \leq \| \int_{t_1}^{t_2} \sum_{j=-1}^{M} \Delta_j \partial_t u(x, t) \, dt \ \|_{L^\infty(\mathbb{R}^n)}$$

$$+ C \sup_{t \in [t_1, t_2]} \sum_{j=M+1}^{\infty} \| \Delta_j \partial_x u(.\,,t) \|_{L^Q(\mathbb{R}^n)} 2^{\frac{jn}{Q} - j}$$

$$\leq \| \int_{t_1}^{t_2} \sum_{j=-1}^{M} \Delta_j \partial_t u(x,t)\, dt \|_{L^Q(\mathbb{R}^n)} 2^{\frac{Mn}{Q}}$$

$$+ C \| \partial_x u(.\,,t) \|_{L^\infty([0,T];\, L^Q(\mathbb{R}^n))} 2^{\frac{Mn}{Q} - M}$$

$$\leq \| \int_{t_1}^{t_2} \partial_t u(.\,,t) \quad dt \|_{L^Q(\mathbb{R}^n)} 2^{\frac{Mn}{Q}}$$

$$+ C \| \partial_x u(.\,,t) \|_{L^\infty([0,T];\, L^Q(\mathbb{R}^n))} 2^{\frac{Mn}{Q} - M}$$

$$\leq C(\int_{t_1}^{t_2} f(t)\, dt)\, 2^{\frac{Mn}{Q}} + C \| \partial_x u(.\,,t) \|_{L^\infty([0,T];\, L^Q(\mathbb{R}^n))} 2^{\frac{Mn}{Q} - M} \,.$$

Choosing $M = [-log_2 \int_{t_1}^{t_2} f(t)\, dt]$ yields (9.17). This concludes the proof of Proposition 9.3. QED.

We now study the equation for the velocity field

$$\begin{cases} \partial_t v_\varepsilon(x,t) = -(v_\varepsilon(x,t), \partial_x) v_\varepsilon(x,t) - \partial_x p_\varepsilon(x,t) + F_\varepsilon(x,t); \ x \in \mathbb{R}^2, \ t \in [0,\infty), \\ curl_x\ v_\varepsilon(x,t) = \omega_\varepsilon(x,t) \,; \ div_x v_\varepsilon(x,t) = 0, \\ v_\varepsilon(x,0) = v_{0\varepsilon}(x). \end{cases} \tag{9.19}$$

Here $v_\varepsilon(x,t)$ is the same as in (9.3), $p_\varepsilon(x,t)$ denotes the pressure;

$$v_{0\varepsilon}(x) = \left( curl^{-1} Re \left( \varepsilon^{\rho + i\zeta + 1 - \alpha} \eta \right) \right)(\varepsilon^{-1} x) + \chi(|x|)\, \varepsilon^{1 - \alpha}\, V(\varepsilon^{-1} x); \tag{9.20}$$

$$F_\varepsilon(x,t) = \varepsilon^{1 - 2\alpha} F(\varepsilon^{-1} x, \varepsilon^{-\alpha} t); \tag{9.21}$$

$$F(x,t) = \partial_t [\mathcal{R}(t)^{1-\alpha} \chi(\varepsilon|x|)\, V(\mathcal{R}(t)^{-1} x)]\,. \tag{9.22}$$

We first compute $F(x,t)$. Using (6.9),

$$F(x,t) = \chi(\varepsilon|x|) \dot{\mathcal{R}}(t) \mathcal{R}(t)^{-\alpha} [(1-\alpha) V(y) - (y, \partial_y) V(y)]|_{y = \mathcal{R}(t)^{-1} x}$$

$$= \gamma\, \chi(\varepsilon|x|) \mathcal{R}(t)^{1-2\alpha} [(1-\alpha) V(y) - (y, \partial_y) V(y)]|_{y = \mathcal{R}(t)^{-1} x}$$

$$= \gamma\, \chi(\varepsilon|x|) \mathcal{R}(t)^{1-2\alpha} [-\alpha \Omega(|y|) y^\perp - |y| \Omega'(|y|) y^\perp]|_{y = \mathcal{R}(t)^{-1} x}$$

$$= \gamma\, \chi(\varepsilon|x|) \mathcal{R}(t)^{-2\alpha} [-\alpha \Omega(|y|) - |y| \Omega'(|y|)]|_{y = \mathcal{R}(t)^{-1} x}\, x^\perp$$

$$= \gamma\, \chi(\varepsilon|x|) \mathcal{R}(t)^{-2\alpha} B(\mathcal{R}(t)^{-1} x) x^\perp \,. \tag{9.23}$$

Here $B(y)$ is the same as in (6.6). From (9.21), (9.23) we derive

$$F_\varepsilon(x,t) = \varepsilon^{1-2\alpha}\,\gamma\,\chi(|x|)(1+\alpha\gamma\varepsilon^{-\alpha}t)^{-2}B\left(\varepsilon^{-1}(1+\alpha\gamma\varepsilon^{-\alpha}t)^{-\frac{1}{\alpha}}x\right)\varepsilon^{-1}x^\perp$$

$$= \gamma\,\chi(|x|)(\varepsilon^\alpha+\alpha\gamma t)^{-2}B((\varepsilon^\alpha+\alpha\gamma t)^{-\frac{1}{\alpha}}x)\,x^\perp\,. \tag{9.24}$$

Theorem 9.3. For any $T>0$ the family of functions $\{v_\varepsilon(x,t)\}$ on $\mathbb{R}^2\times[0,T]$ is equicontinuos as $\varepsilon\in(0,1]$.

Proof. From (9.19), (9.15), (9.24) using boundedness of the Leray projection in $\left(L^Q(\mathbb{R}^2)\right)^2$,

$\forall Q\in(1,\infty)$, we get

$$\|\partial_t v_\varepsilon(.,t)\|_{L^Q(\mathbb{R}^2)}\leq C\|v_\varepsilon(.,t)\|_{L^\infty(\mathbb{R}^2)}\ \|\omega_\varepsilon(.,t)\|_{L^Q(\mathbb{R}^2)}$$

$$+C\|F_\varepsilon(.,t)\|_{L^Q(\mathbb{R}^2)}$$

$$\leq\ C\|\omega_\varepsilon(.,t)\|_{L^{2,1}(\mathbb{R}^2)}\ \|\omega_\varepsilon(.,t)\|_{L^Q(\mathbb{R}^2)}$$

$$+ C\gamma\ \|\chi(|x|)(\varepsilon^\alpha+\alpha\gamma t)^{-2}B((\varepsilon^\alpha+\alpha\gamma t)^{-\frac{1}{\alpha}}x)\,x^\perp\|_{L^Q(\mathbb{R}^2,dx)}$$

$$\leq\ C+C(\varepsilon^\alpha+\alpha\gamma t)^{\frac{1}{\alpha}+\frac{2}{Q\alpha}-2}\ \|B(y)y^\perp\|_{L^Q(\mathbb{R}^2,dy)}$$

with $Q\in(2,\frac{2}{\alpha})$. But $B(y)\,y^\perp\in L^Q(\mathbb{R}^2)$, since $B(y)=C|y|^{-2}$ as $|y|\geq M>0$. Therefore,

$$\|\partial_t v_\varepsilon(.,t)\|_{L^Q(\mathbb{R}^2)}\leq C(1+(\alpha\gamma t)^{\frac{1}{\alpha}+\frac{2}{Q\alpha}-2} \tag{9.25}$$

for $Q\in\left(2,\frac{2}{\alpha}\right),\varepsilon\in(0,1]$ with the constant $C$ uniform with respect to $\varepsilon$.

But the function

$$f(t)=\ C\left(1+(\alpha\gamma t)^{\frac{1}{\alpha}+\frac{2}{Q\alpha}-2}\right)\in L^1_{loc}([0,\infty),dt)\,. \tag{9.26}$$

From the Proposition 9.2 and from the Calderon-Zygmund inequality

$$\partial_x v_\varepsilon(.,t)\in L^\infty_{loc}([0,\infty);L^Q(\mathbb{R}^2)) \tag{9.27}$$

for $Q\in\left(2,\frac{2}{\alpha}\right)$.

Theorem 9.3 now follows from (9.25), (9.26), (9.27) and from the Proposition 9.3, where we select an arbitrary $Q\in\left(2,\frac{2}{\alpha}\right)$. QED

Theorem 9.4. For any $T > 0$ the family of functions $\{v_\varepsilon(x, t)\}$ is uniformly bounded in $L^\infty(\mathbb{R}^2 \times [0, T])$, as $\varepsilon \in (0, 1]$.

Proof. Using Proposition 9.2, Lemma 4.1 and the identity

$$v_\varepsilon(., t) = \mathcal{K} * \omega_\varepsilon(., t),$$

we arrive at the statement of Theorem 9.4. QED.

# § 10. Construction of a non-radial solution

Theorem 10.1 There exists a weak solution $\bar\omega(x, t)$, $\bar v(x, t)$ to the system (4.1) so that

$$\bar\omega(., t) \in L^\infty_{loc}\left([0, T]; L^Q_m(\mathbb{R}^2)\right), \ \ \forall\, Q \in \left(1, \frac{2}{\alpha}\right),$$

$$\bar v(., t) = \mathcal{K}(.) * \bar\omega(., t),$$

$$\bar\omega(x, t) = curl_x \bar v(x, t);$$

$$\bar v(x, t) \in C_m(\mathbb{R}^2 \times [0, T]) \cap L^\infty(\mathbb{R}^2 \times [0, T]) \ for \ any \ T > 0,$$

$\bar v(x, t)$ is Holder continuous on $\mathbb{R}^2 \times [0, T]$ for any $T > 0$ with an arbitrary Holder exponent $s \in (0, 1 - \alpha)$. This weak solution satisfies the integral identity

$$\int \overline{\omega_0}(x)\, \psi(x, 0)\, dx dt + \int\limits_0^\infty \int \bar\omega(x, t)\, \partial_t\, \psi(x, t)\, dx\, dt$$

$$+ \int\limits_0^\infty \int \bar\omega(x, t)\, (\bar v(x, t), \partial_x)\psi(x, t)\, dx\, dt + \int\limits_0^\infty \int \bar Z(x, t)\, \psi(x, t)\, dx\, dt \ = 0$$

valid for any test function $\psi(x, t) \in C^\infty\left([0, \infty); S'(\mathbb{R}^2)\right)$. The external force $\bar Z(x, t)$ is given by

$$\bar Z(x, t) = \gamma(\alpha\gamma t)^{-2}(\chi(|x|)W((\alpha\gamma t)^{-\frac{1}{\alpha}}\, x) + |x|\chi'(|x|)B((\alpha\gamma t)^{-\frac{1}{\alpha}}\, x))$$

$$= \frac{\partial}{\partial t} curl_x\left[\chi(|x|)(\alpha\gamma t)^{-1+\frac{1}{\alpha}}\, V\left((\alpha\gamma t)^{-\frac{1}{\alpha}}\, x\right)\right],$$

$$\overline{\omega_0}(x) = \frac{1}{2-\alpha} \, curl \, [(\chi(|x|) \, |x|^{-\alpha} x^\perp].$$

Proof. We will pass to the limit as $\varepsilon \to 0$ in the integral identity

$$\int \omega_{0\varepsilon}(x) \, \psi(x,0) \, dx \qquad\qquad (I_1)$$

$$+ \int\limits_0^\infty \int \omega_\varepsilon(x,t) \, \partial_t \psi(x,t) \, dx \, dt \qquad\qquad (I_2)$$

$$+ \int\limits_0^\infty \int \omega_\varepsilon(x,t) \, (v_\varepsilon(x,t), \partial_x) \psi(x,t) \, dx \, dt \qquad\qquad (I_3)$$

$$+ \int_0^\infty \int Z_\varepsilon(x,t)$$
$$\psi(x,t) \, dx \, dt \qquad\qquad (I_4)$$

$$= 0 \qquad\qquad (10.1)$$

We start with the integral $I_1(\varepsilon)$.

From (9.5)

$$\omega_{0\varepsilon}(x) = Re(\varepsilon^{\rho+i\zeta-\alpha} \eta(\varepsilon^{-1}x)) + \, curl_x \, [\, \chi(|x|) \, \varepsilon^{1-\alpha} \, V(\varepsilon^{-1}x) \,]$$

$$= Re \, (\varepsilon^{\rho+i\zeta-\alpha} \eta(\varepsilon^{-1}x)) + \varepsilon^{-\alpha} \chi(|x|) \, (curl_x V)(\varepsilon^{-1}x)$$

$$+ |x|\chi'(|x|) \, \varepsilon^{-\alpha} \, \Omega(\varepsilon^{-1}x).$$

Since $\rho > 0$, $Re(\varepsilon^{\rho+i\zeta-\alpha} \eta(\varepsilon^{-1}x)) \to 0$ strongly in any $L_m^Q(\mathbb{R}^2)$ with $Q \in \left(2, \frac{2}{\alpha}\right)$. Therefore,

$$\omega_{0\varepsilon}(x) \to \frac{1}{2-\alpha} curl \, [(\chi(|x|) \, |x|^{-\alpha} x^\perp] \qquad\qquad (10.2)$$

in $S_m'(\mathbb{R}^2)$ as $\varepsilon \to 0+$.

From (10.1), (10.2)

$$\lim_{\varepsilon \to 0+} I_1(\varepsilon) = \int \overline{\omega_0}(x) \, \psi(x,0) \, dx, \qquad\qquad (10.3)$$

where

$$\overline{\omega_0}(x) = \frac{1}{2-\alpha} curl \, [(\chi(|x|) \, |x|^{-\alpha} x^\perp]. \qquad\qquad (10.4)$$

We turn to the second term $I_2(\varepsilon)$:

$$I_2(\varepsilon) = \int_0^\infty \int \omega_\varepsilon(x,t) \; \partial_t \psi(x,t) \; dx \, dt \, . \tag{10.5}$$

Using the Proposition 9.2, we select a sequence $\{\varepsilon_j\}_{j=1}^\infty \quad \varepsilon_j \in (0,1] \;\; \forall j = 1,2,3,\dots; \; \varepsilon_j \to 0$ as $j \to \infty$, so that

$$\omega_{\varepsilon_j}(.,t) \to \bar\omega(.,t) \; weak -* \; in \; L^\infty\left([0,T]; \; L_m^Q(\mathbb{R}^2)\right) \tag{10.6}$$

for any $T > 0$ and for any $Q \in \left(1, \frac{2}{\alpha}\right)$. Thus the limit

$$\bar\omega(.,t) \in L_{loc}^\infty\left([0,T]; \; L_m^Q(\mathbb{R}^2)\right), \;\; \forall \, Q \in \left(1, \frac{2}{\alpha}\right). \tag{10.7}$$

Therefore,

$$\lim_{j\to\infty} I_2(\varepsilon_j) = \int_0^\infty \int \bar\omega(x,t) \; \partial_t \psi(x,t) \; dx \, dt. \tag{10.8}$$

We now turn to the integral

$$I_3(\varepsilon) = \int\limits_0^\infty \int \omega_\varepsilon(x,t) \; (v_\varepsilon(x,t), \partial_x)\psi(x,t) \; dx \, dt \; . \tag{10.9}$$

From the previous construction by selecting a subsequence of $\{\varepsilon_j\}_{j=1}^\infty$, we may and will assume using the Arzela-Ascoli theorem and theorems 9.3 and 9.4 that

$$v_{\varepsilon_j}(x,t) \to \bar v(x,t) \; uniformly \; on \; every \; cillinder \; B_N \times [0,T] \; , \tag{10.10}$$

for any $N = 1,2,3,\dots, \;\; B_N = \{x \in \mathbb{R}^2 | \; |x| \le N\}$, and for any $T > 0$. We also have

$$\bar\omega(x,t) = curl_x \bar v(x,t); \; \bar v(x,t) \in C_m(\mathbb{R}^2 \times [0,T]) \cap L^\infty(\mathbb{R}^2 \times [0,T]),$$

Holder continuous with the Holder exponent $1 - \frac{2}{Q}$ , $\forall Q \in \left(2, \frac{2}{\alpha}\right)$, on $\mathbb{R}^2 \times [0,T]$ , for any $T > 0$.

Given this information, we can pass to the limit

$$\lim_{j\to\infty} I_3(\varepsilon_j) = \int\limits_0^\infty \int \bar\omega(x,t) \; (\bar v(x,t), \partial_x)\psi(x,t) \; dx \, dt \; . \tag{10.11}$$

It remains to study the limit of the integral

$$I_4(\varepsilon_j) = \int_0^\infty \int Z_{\varepsilon_j}(x,t) \; \psi(x,t) \; dx \, dt.$$

Using Proposition 9.1 we get

$$\lim_{j\to\infty} I_3(\varepsilon_j) = \int_0^\infty \int \bar{Z}(x,t)\,\psi(x,t)\,dx\,dt, \tag{10.12}$$

where

$$\bar{Z}(x,t) = \gamma(\alpha\gamma t)^{-2}(\chi(|x|)W((\alpha\gamma t)^{-\frac{1}{\alpha}}\,x) + |x|\chi'(|x|)B((\alpha\gamma t)^{-\frac{1}{\alpha}}\,x))$$

$$= \frac{\partial}{\partial t}\,curl_x\left[\chi(|x|)(\alpha\gamma t)^{-1+\frac{1}{\alpha}}\,V\left((\alpha\gamma t)^{-\frac{1}{\alpha}}\,x\right)\right]. \tag{10.13}$$

From (10.1), (10.3), (10.8), (10.11), (10.12) by passing to the limit as $j \to \infty$ we obtain the following

Identity:

$$\int \overline{\omega_0}(x)\,\psi(x,0)\,dx + \int_0^\infty \int \bar{\omega}(x,t)\,\partial_t\,\psi(x,t)\,dx\,dt$$

$$+ \int_0^\infty \int \bar{\omega}(x,t)\,(\bar{v}(x,t),\partial_x)\psi(x,t)\,dx\,dt + \int_0^\infty \int \bar{Z}(x,t)\,\psi(x,t)\,dx\,dt = 0 \tag{10.14}$$

valid for any test function $\psi(x,t) \in C^\infty\big([0,\infty); S'(\mathbb{R}^2)\big)$.

This concludes the proof of Theorem 10.1. QED.

Theorem 10.2. The solution $(\bar{\omega}(x,t), \bar{v}(x,t))$ of the Euler equations is not rotationally invariant.

Proof. We study in more detail the limit (10.6).

From (6.3), (9.1), (8.1), (8.2)

$$\omega_\varepsilon(x,t) = (\varepsilon^\alpha + \alpha\gamma t)^{-1}\sigma_1\left(y, \frac{1}{\alpha\gamma}\log(1+\alpha\gamma\varepsilon^{-\alpha}t)\right)\big|_{y=\frac{x}{(\varepsilon^\alpha+\alpha\gamma t)^{\frac{1}{\alpha}}}}$$

$$+(\varepsilon^\alpha + \alpha\gamma t)^{-1}[\chi(|x|)\,(curl_y V)\left(\frac{x}{(\varepsilon^\alpha+\alpha\gamma t)^{\frac{1}{\alpha}}}\right) + |x|\chi'(|x|)\Omega(\frac{x}{(\varepsilon^\alpha+\alpha\gamma t)^{\frac{1}{\alpha}}})] \tag{10.15}$$

$$= J_\varepsilon(x,t) + K_\varepsilon(x,t)\quad.$$

Then,

$$K_\varepsilon(x,t) \to (\alpha\gamma t)^{-1}[\chi(|x|)\,(curl_y V)\left(\frac{x}{(\alpha\gamma t)^{\frac{1}{\alpha}}}\right) + |x|\chi'(|x|)\Omega(\frac{x}{(\alpha\gamma t)^{\frac{1}{\alpha}}})]$$

$$= (\alpha\gamma t)^{-1}\,curl_y\left[\chi\left((\alpha\gamma t)^{\frac{1}{\alpha}}|y|\right)V(|y|)\right]\big|_{y=\frac{x}{(\alpha\gamma t)^{\frac{1}{\alpha}}}},$$

as $\varepsilon \to 0$, strongly in $L^\infty([0,T]; L_m^Q(\mathbb{R}^2))$ for any $T > 0$ and any $Q \in \left(1, \frac{2}{\alpha}\right)$.

To study the limit as $\varepsilon \to 0$ of the first term in the right side of (10.15), we use the theorem 8.2 and (8.34)

$$\|\sigma_1(\tau) - e^{\tau L}\sigma_0\|_{L^Q(\mathbb{R}^2)} \leq CC_0^{\frac{\alpha}{\rho}}\|e^{\tau L}\sigma_0\|_{L^Q(\mathbb{R}^2)}, \qquad \forall \tau \in [0, T_*], \tag{10.16}$$

where

$$T_* = \frac{1}{a}\left(\rho \log \varepsilon^{-1} - log \frac{C_1}{C_0}\right).$$

Notice that

$$(e^{\tau L}\sigma_0)(y, t) = Re\left(\varepsilon^{\rho + i\zeta} e^{\lambda \tau(t)}\eta(y)\right).$$

We now set

$$\rho + i\zeta = \gamma^{-1}\lambda \ . \tag{10.17}$$

Therefore,

$$\|\varepsilon^{-\alpha}\mathcal{R}(\varepsilon^{-\alpha}t)^{-\alpha}\sigma_1(\varepsilon^{-1}\mathcal{R}(\varepsilon^{-\alpha}t)^{-1}x, \tau(\varepsilon^{-\alpha}t))$$

$$- (\varepsilon^\alpha + \alpha\gamma t)^{-1} Re\left(\varepsilon^{\rho + i\zeta} e^{\lambda\frac{1}{\alpha\gamma}\log(1 + \alpha\gamma\varepsilon^{-\alpha}t)}\eta\left(\frac{x}{(\varepsilon^\alpha + \alpha\gamma t)^{\frac{1}{\alpha}}}\right)\right)\|_{L^Q(\mathbb{R}^2)}$$

$$\leq CC_0^{\frac{\alpha}{\rho}}\|(\varepsilon^\alpha + \alpha\gamma t)^{-1} Re\left(\varepsilon^{\rho + i\zeta} e^{\lambda\frac{1}{\alpha\gamma}\log(1 + \alpha\gamma\varepsilon^{-\alpha}t)}\eta\left(\frac{x}{(\varepsilon^\alpha + \alpha\gamma t)^{\frac{1}{\alpha}}}\right)\right)\|_{L^Q(\mathbb{R}^2)}, \tag{10.18}$$

as long as

$$\frac{1}{\alpha\gamma}\log(1 + \alpha\gamma\varepsilon^{-\alpha}t) \in [0, T_*] = [0.\frac{1}{a}\left(\rho \log \varepsilon^{-1} - log \frac{C_1}{C_0}\right). \tag{10.19}$$

But

$$(\varepsilon^\alpha + \alpha\gamma t)^{-1} Re\left(\varepsilon^{\rho + i\zeta} e^{\lambda\frac{1}{\alpha\gamma}\log(1 + \alpha\gamma\varepsilon^{-\alpha}t)}\eta\left(\frac{x}{(\varepsilon^\alpha + \alpha\gamma t)^{\frac{1}{\alpha}}}\right)\right)$$

$$= (\varepsilon^\alpha + \alpha\gamma t)^{-1} Re\left((\varepsilon^\alpha + \alpha\gamma t)^{\frac{\lambda}{\alpha\gamma}}\eta\left(\frac{x}{(\varepsilon^\alpha + \alpha\gamma t)^{\frac{1}{\alpha}}}\right)\right) \ . \tag{10.20}$$

From (10.18) --(10.20)

$$\|\varepsilon^{-\alpha}\mathcal{R}(\varepsilon^{-\alpha}t)^{-\alpha}\sigma_1(\varepsilon^{-1}\mathcal{R}(\varepsilon^{-\alpha}t)^{-1}x, \tau(\varepsilon^{-\alpha}t))$$

$$- (\varepsilon^{\alpha} + \alpha\gamma t)^{-1} Re\left((\varepsilon^{\alpha} + \alpha\gamma t)^{\frac{\lambda}{\alpha\gamma}}\eta\left(\frac{x}{(\varepsilon^{\alpha} + \alpha\gamma t)^{\frac{1}{\alpha}}}\right)\right)\|_{L^Q(\mathbb{R}^2)}$$

$$\leq C C_0^{\frac{\alpha}{\rho}} C_1^{1-\frac{\alpha}{\rho}} \|(\varepsilon^{\alpha} + \alpha\gamma t)^{-1} Re\left((\varepsilon^{\alpha} + \alpha\gamma t)^{\frac{\lambda}{\alpha\gamma}}\eta\left(\frac{x}{(\varepsilon^{\alpha} + \alpha\gamma t)^{\frac{1}{\alpha}}}\right)\right)\|_{L^Q(\mathbb{R}^2)}. \tag{10.21}$$

Since $Re\lambda > 0$, and since $\eta(y) \in L_m^1(\mathbb{R}^2) \cap L_m^\infty(\mathbb{R}^2)$,

$$(\varepsilon^{\alpha} + \alpha\gamma t)^{-1} Re\left((\varepsilon^{\alpha} + \alpha\gamma t)^{\frac{\lambda}{\alpha\gamma}}\eta\left(\frac{x}{(\varepsilon^{\alpha} + \alpha\gamma t)^{\frac{1}{\alpha}}}\right)\right) \to (\alpha\gamma t)^{-1} Re\left((\alpha\gamma t)^{\frac{\lambda}{\alpha\gamma}}\eta\left(\frac{x}{(\alpha\gamma t)^{\frac{1}{\alpha}}}\right)\right)$$

as $\varepsilon \to 0$ strongly in any space $L^\infty([0,T]; L_m^Q(\mathbb{R}^2))$, $\forall T > 0$, $\forall Q \in \left(2, \frac{2}{\alpha}\right)$. Therefore,

$$w-* \lim_{j\to\infty} J_{\varepsilon_j}(x,t) = (\alpha\gamma t)^{-1} Re\left((\alpha\gamma t)^{\frac{\lambda}{\alpha\gamma}}\eta\left(\frac{x}{(\alpha\gamma t)^{\frac{1}{\alpha}}}\right)\right) + \bar{H}(x,t) \tag{10.22}$$

in any space $L^\infty([0,T]; L_m^Q(\mathbb{R}^2))$, $\forall Q \in \left(2, \frac{2}{\alpha}\right)$ provided (see (10.19))

$$\frac{1}{\alpha\gamma}\log(\alpha\gamma T) < -\log\frac{C_1}{C_0}. \tag{10.23}$$

Here $\bar{H}(x,t)$ satisfies

$$\|\bar{H}(x,t)\|_{L^\infty([0,T]; L_m^Q(\mathbb{R}^2))}$$

$$\leq C C_0^{\frac{\alpha}{\rho}} C_1^{1-\frac{\alpha}{\rho}} \|(\alpha\gamma t)^{-1} Re\left((\alpha\gamma t)^{\frac{\lambda}{\alpha\gamma}}\eta\left(\frac{x}{(\alpha\gamma t)^{\frac{1}{\alpha}}}\right)\right)\|_{L^\infty([0,T]; L_m^Q(\mathbb{R}^2))}, \tag{10.24}$$

with $T$ satisfing (10.23) and $\forall Q \in \left(2, \frac{2}{\alpha}\right)$. Selecting $C_0 > 0$ sufficiently small, we conclude that

$$w-* \lim_{j\to\infty} J_{\varepsilon_j}(x,t)$$

in any $L^\infty([0,T]; L^Q_m(\mathbb{R}^2))$ space is not rotation invariant for $T > 0$ , $\forall Q \in \left(2, \frac{2}{\alpha}\right)$. This concludes the proof of Theorem 8.2. QED.

Remark. There is a solution with the same initial data

$$\overline{\omega_0}(x) = \frac{1}{2-\alpha} curl \left[\left(\chi(|x|) |x|^{-\alpha} x^\perp\right].$$

and with same

$$\bar{Z}(x,t) = \frac{\partial}{\partial t} curl_x \left[\chi(|x|)(\alpha\gamma t)^{-1+\frac{1}{\alpha}} V\left((\alpha\gamma t)^{-\frac{1}{\alpha}} x\right)\right],$$

given by

$$\overline{\omega}(x,t) = (\alpha\gamma t)^{-1} curl_y \left[\chi\left((\alpha\gamma t)^{\frac{1}{\alpha}}|y|\right) V(|y|)\right]\Big|_{y=\frac{x}{(\alpha\gamma t)^{\frac{1}{\alpha}}}} ;$$

$$\bar{v}(x,t) = (\alpha\gamma t)^{-1+\frac{1}{\alpha}} \chi(|x|) V\left((\alpha\gamma t)^{-\frac{1}{\alpha}} x\right).$$

Therefore, uniqueness fails.